\documentclass[11pt]{article}
\usepackage[left=1in,right=1in,top=1in,bottom=1in]{geometry}
\usepackage{times}
\usepackage{expl3}
\usepackage{cite}
\usepackage[table]{xcolor}
\usepackage{multirow}
\usepackage{stackengine} 
\usepackage{hhline}
\usepackage{lipsum}
\usepackage{titlesec}
\usepackage{wrapfig}
\usepackage{epsfig}
\usepackage{graphicx}
\usepackage{amsmath}
\usepackage{amssymb}
\usepackage{epstopdf}
\usepackage{boldline}
\usepackage{calligra}
\usepackage{url}
\usepackage{blindtext}

\newcommand{\define}{\stackrel{\mbox{\tiny def}}{=}}

\newtheorem{definition}{Definition}
\newtheorem{theorem}{Theorem}

\newtheorem{lemma}{Lemma}

\usepackage{mathtools}
\usepackage{epstopdf}
\usepackage{balance}
\usepackage{thmtools}
\usepackage{thm-restate}
\usepackage{hyperref}
\usepackage{cleveref}

\usepackage[ruled,vlined]{algorithm2e}
\include{pythonlisting}
\newcommand{\ostar}{\mathbin{\mathpalette\make@circled\star}}

\makeatletter
\newcommand{\removelatexerror}{\let\@latex@error\@gobble}
\makeatother
\makeatletter
\newcommand*{\rom}[1]{\expandafter\@slowromancap\romannumeral #1@}
\makeatother

\ExplSyntaxOn
\newcommand\latinabbrev[1]{
  \peek_meaning:NTF . {
    #1\@}%
  { \peek_catcode:NTF a {
      #1.\@ }%
    {#1.\@}}}
\ExplSyntaxOff

\titleclass{\subsubsubsection}{straight}[\subsubsection]

\begin{document}
\vspace{1cm}
\title{Generalized  T-product Tensor Bernstein Bounds}\vspace{1.8cm}
\author{Shih~Yu~Chang 
\thanks{Shih Yu Chang is with the Department of Applied Data Science,
San Jose State University, San Jose, CA, U. S. A. (e-mail: {\tt
shihyu.chang@sjsu.edu}).
           }
Yimin Wei
\thanks{1This work was supported by the National Natural Science Foundation of China (No. 11771099) and
Innovation Program of Shanghai Municipal Education Commission.}
\thanks{Yimin Wei is with the School of Mathematical Sciences,
Shanghai Key Laboratory of Contemporary Applied Mathematics, 
Fudan University, 
Shanghai, 200433, PR China(e-mail: {\tt
ymwei@fudan.edu.cn}).
           }}
\maketitle


\begin{abstract}
Since Kilmer et al. introduced the new multiplication method between two third-order tensors around 2008 and third-order tensors with such multiplication structure are also called as T-product tensors, T-product tensors have been applied to many fields in science and engineering, such as low-rank tensor approximation, signal processing, image feature extraction, machine learning, computer vision, and the multi-view clustering problem, etc. However, there are very few works dedicated to exploring the behavior of random T-product tensors. This work considers the problem about the tail behavior of the unitarily invariant norm for the summation of random symmetric T-product tensors. Majorization and antisymmetric Kronecker product tools are main techniques utilized to establish inequalities for unitarily norms of multivariate T-product tensors. The Laplace transform method is integrated with these inequalities for unitarily norms of multivariate T-product tensors to provide us with Bernstein Bounds estimation of Ky Fan $k$-norm for functions of the symmetric random T-product tensors summation. 
\end{abstract}


\section{Introduction}\label{sec:Introduction} 


Since Kilmer et al. introduced the new multiplication method between two third-order tensors (T-product tensors), many new algebraic properties about such new multiplication rule between two third-order tensors are investigated recently~\cite{kilmer2008third, kilmer2013third}. For example, the singular value decomposition (SVD) for third-order tensors via the tensor T-product is proposed in~\cite{kilmer2011factorization}. Some authors suggest a new framework by treating third-order tensors as linear operators on a space of matrices, see ~\cite{braman2010third}. In~\cite{miao2021t}, many useful tools of linear algebra are extended to the third-order tensors, including the T-Jordan canonical form, tensor decomposition theory, T-group inverse and T-Drazin inverse, and so on. Moreover, the authors in~\cite{lund2020tensor} proposed a definition of tensor functions based on the T-product of third-order F-square tensors, and Miao, Qi and Wei generalized the tensor T-function from F-square third-order tensors to rectangular tensors in~\cite{miao2020generalized}. These useful algebraic properties of T-product tensors have been discovered as powerful tools in many science and engineering fields: signal processing~\cite{zhang2016exact, semerci2014tensor}, machine learning~\cite{settles2007multiple}, computer vision~\cite{zhang2014novel, martin2013order}, image processing~\cite{khalil2021efficient}, low-rank tensor approximation~\cite{xu2013parallel, zhou2017tensor, qi2021tsingular}, etc. 

Although T-product tensors have attracted many practical applications, all of these applications of T-product tensors assume that T-product tensors under consideration are deterministic. This assumption is not practical in general scientific and engineering applications based on T-product tensors. In~\cite{chang2021TProdI,chang2021TProdII}, the authors have tried to establish several new tail bounds for sums of random T-product tensors. These probability bounds characterize large-deviation behavior of the extreme T-eigenvalue of the sums of random T-product tensors~\footnote{Definitions about T-eigenvalues and T-singular values associated to T-product tensors are given in Section ~\ref{sec:T-product Tensor Fundamental Facts}.}. The authors first apply Lapalace transform method and Lieb's concavity theorem for T-product tensors obtained from the work~\cite{chang2021TProdI} to build several inequalities based on random T-product tensors, then utilize these inequalities to generalize the classical bounds associated with the names Chernoff, and Bernstein from the scalar to the T-product tensor setting. Tail bounds for the norm of a sum of random rectangular T-product tensors are also derived from corollaries of random symmetric T-product tensors cases. The proof mechanism is also applied to T-product tensor valued martingales and T-product tensor-based Azuma, Hoeffding and McDiarmid inequalities are also derived~\cite{chang2021TProdII}. 

In this work, we will apply majorization techniques to establish new Bernstein bounds based on the summation of random symmetric T-product tensors. Compared to the previous work studied in~\cite{chang2021TProdI,chang2021TProdII}, we make following generalizations: (1) besides bounds related to extreme values of T-eigenvalues, we consider more general unitarily invariant norm for T-product tensors; (2) the bounds derived in ~\cite{chang2021TProdII} can only be applied to the identity map for the summation of random symmetric T-product tensors, this work can derive new bounds for any polynomial function raised by any power greater or equal than one for the summation of random symmetric T-product tensors. In order to drive these new bounds, we also establish Courant-Fischer min-max theorem for T-product tensors in Theorem~\ref{thm:Courant-Fischer T-product} and marjoization relation for T-singular values in Lemma~\ref{lma:singular values major relation}. Our main theorem is provided below:
\begin{restatable}[Generalized T-product Tensor Bernstein Bound]{thm}{thmGeneralizedTensorBernsteinBound}\label{thm:Generalized Tensor Bernstein Bound}
Consider a sequence $\{ \mathcal{X}_j  \in \mathbb{R}^{m \times m \times p} \}$ of independent, random symmetric T-product tensors with random structure defined by Definition~\ref{def:Random structure of discussed random symmetric T-product tensor}. Let $g$ be a polynomial function with degree $n$ and nonnegative coefficients $a_0, a_1, \cdots, a_n$ raised by power $s \geq 1$, i.e., $g(x) = \left(a_0 + a_1 x  +\cdots + a_n x^n \right)^s$ with $s \geq 1$. Suppose following condition is satisfied:
\begin{eqnarray}\label{eq:special cond Bernstein Bound}
g \left( \exp\left(t \sum\limits_{j=1}^{m} \mathcal{X}_j \right)\right)  \succeq \exp\left(t g \left( \sum\limits_{j=1}^{m} \mathcal{X}_j   \right) \right)~~\mbox{almost surely},
\end{eqnarray}
where $t > 0$, and we also have 
\begin{eqnarray}
\mathcal{X}^p_j \preceq \frac{p! \mathcal{A}^2}{2}
\mbox{~~ almost surely for $p=2,3,4,\cdots$.}
\end{eqnarray}
Then we have following inequality:
\begin{eqnarray}\label{eq1:thm:Generalized Tensor Bernstein Bound}
\mathrm{Pr} \left( \left\Vert g\left( \sum\limits_{j=1}^{m} \mathcal{X}_j  \right)\right\Vert_{(k)}  \geq \theta \right) & \leq & (n+1)^{s-1} \inf\limits_{t > 0} e^{- \theta t }k \cdot \nonumber \\
&  & \left\{a_0^s + \sum\limits_{l=1}^{n} a_l^{l s}\left[1 + mlst \Phi(m, d_1, d_2) + 
\frac{(mlst)^2   \sigma_1(\mathcal{A}^2) }{2(1 - mlst)  }
   \right] \right\}.
\end{eqnarray}
\end{restatable}

The rest of this paper is organized as follows. In Section~\ref{sec:T-product Tensors} , we review T-product tensors basic concepts and introduce a powerful scheme about antisymmetric  Kronecker product for T-product tensors. In Section~\ref{sec:Multivariate T-product Tensor Norm Inequalities}, we apply a majorization technique to prove T-product tensor norm inequalities. We then apply new derived T-product tensor norm inequalities to obtain random T-product tensor Bernstein bounds for the extreme T-eigenvalues and Ky Fan $k$-norm in Section~\ref{sec:Applications of T-product Tensor Norm Inequalities}. Finally, concluding remarks are given by Section~\ref{sec:Conclusions}.

\section{T-product Tensors}\label{sec:T-product Tensors} 

In this section, we will introduce fundamental facts about T-product tensors in Section~\ref{sec:T-product Tensor Fundamental Facts}. Several unitarily invariant norms about a T-product tensor are defined in Section~\ref{sec:Unitarily Invariant T-product Tensor Norms}. A powerful scheme about antisymmetric Kronecker product for T-product tensors will be provided by Section~\ref{sec:Antisymmetric  Kronecker Product for T-product Tensors}.

\subsection{T-product Tensor Fundamental Facts}\label{sec:T-product Tensor Fundamental Facts}

For a third order tensor $\mathcal{C} \in \mathbb{R}^{m \times n \times p}$, we define $\mbox{bcirc}$ operation to the tensor $\mathcal{C}$ as:
\begin{eqnarray}
\mbox{bcirc} (\mathcal{C} ) \define \left[
    \begin{array}{ccccc}
       \mathbf{C}^{(1)}  &  \mathbf{C}^{(p)}  &  \mathbf{C}^{(p-1)}  & \cdots  &  \mathbf{C}^{(2)}   \\
       \mathbf{C}^{(2)}  &  \mathbf{C}^{(1)}  &  \mathbf{C}^{(p)}  & \cdots  &  \mathbf{C}^{(3)}   \\
       \vdots  &  \vdots  & \vdots  & \cdots  & \vdots   \\
       \mathbf{C}^{(p)}  &  \mathbf{C}^{(p-1)}  &  \mathbf{C}^{(p-2)}  & \cdots  &  \mathbf{C}^{(1)}   \\
    \end{array}
\right],
\end{eqnarray}
where $\mathbf{C}^{(1)}, \cdots, \mathbf{C}^{(p)} \in \mathbb{C}^{ m \times n}$ are frontal slices of tensor $\mathcal{C}$. The inverse operation of $\mbox{bcirc}$ is denoted as $\mbox{bcirc}^{-1}$ with relation $\mbox{bcirc}^{-1} ( \mbox{bcirc} ( \mathcal{C} )) \define \mathcal{C}$. Another operation to the tensor $\mathcal{C}$ is unfolding, denoted as $\mbox{unfold} (\mathcal{C} )$, which is defined as: 
\begin{eqnarray}
\mbox{unfold} (\mathcal{C} ) \define \left[
    \begin{array}{c}
       \mathbf{C}^{(1)}  \\
       \mathbf{C}^{(2)}  \\
       \vdots    \\
       \mathbf{C}^{(p)}   \\
    \end{array}
\right].
\end{eqnarray}
The inverse operation of $\mbox{unfold}$ is denoted as $\mbox{fold}$ with relation $\mbox{fold}( \mbox{unfold} ( \mathcal{C} )) \define \mathcal{C}$.

The multiplication between two third order tensors, $\mathcal{C} \mathbb{R}^{m \times n \times p}$ and $\mathcal{D} \mathbb{R}^{n \times l \times p} $, is via \emph{T-product} and this multiplication is defined as:
\begin{eqnarray}\label{eq:T prod def}
\mathcal{C} \star \mathcal{D} \define \mbox{fold}\left( \mbox{bcirc}(\mathcal{C}) \cdot \mbox{unfold}(\mathcal{C})\right),
\end{eqnarray}
where $\cdot$ is the standard matrix multiplication. For given third order tensors, if we apply \emph{T-product} to multiply them, we call them \emph{T-product tensors}. A T-product tensor $\mathcal{C} \in \mathbb{R}^{m \times n \times p}$ will be named as square T-product tensor if $m=n$. 

For a symmetric T-product tensor $\mathcal{C} \in \mathbb{R}^{m \times m \times p}$, we define Hermitian transpose of $\mathcal{C}$, denoted by $\mathcal{C}^{\mathrm{H}}$ , as 
\begin{eqnarray}\label{eq:Hermitian Transpose Def}
\mathcal{C}^{\mathrm{H}} = \mbox{bcirc}^{-1}(     (\mbox{bcirc}(\mathcal{C}))^{\mathrm{H}}  ). 
\end{eqnarray}
And a tensor $\mathcal{D} \in \mathbb{C}^{m \times m \times p }$ is called a Hermitian T-product tensor if $ \mathcal{D}^{\mathrm{H}}  = \mathcal{D}$. Similarly, for a symmetric T-product tensor $\mathcal{C} \in \mathbb{R}^{m \times m \times p}$, we define \emph{transpose} of $\mathcal{C}$, denoted by $\mathcal{C}^{T}$ , as 
\begin{eqnarray}\label{eq:Transpose Def}
\mathcal{C}^{T} = \mbox{bcirc}^{-1}(     (\mbox{bcirc}(\mathcal{C}))^{T}  ). 
\end{eqnarray}
And a tensor $\mathcal{D} \in \mathbb{R}^{m \times m \times p }$ is called a symmetric T-product tensor if $ \mathcal{D}^{T}  = \mathcal{D}$. 

The identity tensor $\mathcal{I}_{m,m,p} \in \mathbb{R}^{m \times m \times p }$ can be defined as:
\begin{eqnarray}\label{eq:I_mmp def}
\mathcal{I}_{m,m,p} = \mbox{bcirc}^{-1}(  \mathbf{I}_{mp}  ),
\end{eqnarray}
where $\mathbf{I}_{mp} $ is the identity matrix in $\mathbb{R}^{mp \times mp}$. For a square T-product tensor, $
\mathcal{C} \in \mathbb{R}^{m \times m \times p }$, we say that $\mathcal{C}$ is nonsingular if it has an inverse tensor $\mathcal{D} \in \mathbb{R}^{m \times m \times p } $ such that 
\begin{eqnarray}\label{eq:inverse tensor def}
\mathcal{C} \star \mathcal{D} = \mathcal{D} \star \mathcal{C} = \mathcal{I}_{m, m, p}.
\end{eqnarray}
A zero tensor, denoted as $\mathcal{O}_{mnp} \in \mathbb{C}^{m \times n \times p}$, is a tensor that all elements inside the tensor as $0$. 

For any circular matrix $\mathbf{C} \in \mathbb{R}^{m \times m}$, it can be diagonalized with the normalized Discrete Fourier Transform (DFT) marix, i.e., $\mathbf{C} = \mathbf{F}^{\mathrm{H}}_m \mathbf{D} \mathbf{F}_m   $, where $\mathbf{F}_m$ is the Fourier matrix of size $m \times m$ defined as 
\begin{eqnarray}\label{eq:DFT def}
\mathbf{F}_m \define  \frac{1}{\sqrt{m}} \left[
    \begin{array}{ccccc}
       1 &  1 & 1 & \cdots  & 1  \\
       1 &  \omega  &  \omega^2  & \cdots  & \omega^{(m-1)}  \\
       \vdots  &  \vdots  & \vdots  & \cdots  & \vdots   \\
       1 &   \omega^{(m-1)} &  \omega^{2(m-1)} & \cdots  &   \omega^{(m-1)(m-1)}   \\
    \end{array}
\right],
\end{eqnarray}
where $\omega = \exp(\frac{2 \pi \iota}{m})$ with $\iota^2 = -1$. This DFT matrix can also be used to diagonalize a T-product tensor as~\cite{kilmer2013third}
\begin{eqnarray}\label{eq:block diagonalized format}
\mbox{bcirc}(\mathcal{C}) = \left( \mathbf{F}^{\mathrm{H}}_m \otimes \mathbf{I}_m \right) \mbox{Diag}\left( \mathbf{C}_i: i \in \{1, \cdots, m \} \right)  \left( \mathbf{F}_m \otimes \mathbf{I}_m \right),
\end{eqnarray}
where $ \otimes$ is  the Kronecker Product and $ \mbox{Diag}\left( \mathbf{C}_i: i \in \{1, \cdots, m \} \right) \in \mathbb{C}^{mp \times mp}$ is a diagonal block matrix with the $i$-th diagonal block as the matrix $\mathbf{A}_i$. 

The inner product between two T-product tensors $\mathcal{C} \in \mathbb{C}^{m \times n \times p}$ and $\mathcal{D} \in \mathbb{C}^{m \times n \times p}$ is defined as:
\begin{eqnarray}
\langle \mathcal{C}, \mathcal{D} \rangle = \sum\limits_{i, j, k} c^{\ast}_{i, j, k}d_{i, j, k},
\end{eqnarray}
where $\ast$ is the complex conjugate operation.

We say that a symmetric T-product tensor $\mathcal{C} \in \mathbb{R}^{m \times m \times p}$ is a T-positive definite (TPD) tensor if we have 
\begin{eqnarray}\label{eq: TPD def}
\langle \mathcal{X}, \mathcal{C} \star \mathcal{X} \rangle > 0,
\end{eqnarray}
holds for any non-zero T-product tensor $\mathcal{X} \in \mathbb{R}^{m \times 1 \times p}$. Also, we said that a symmetric T-product tensor is a T-positive semidefinite (TPSD) tensor if we have 
\begin{eqnarray}\label{eq: TPSD def}
\langle \mathcal{X}, \mathcal{C} \mathcal{X} \rangle \geq 0,
\end{eqnarray}
holds for any non-zero T-product tensor $\mathcal{X} \in \mathbb{R}^{m \times 1 \times p}$. Given two T-product tensors $\mathcal{C}, \mathcal{D}$, we use $\mathcal{C} \succ (\succeq) \mathcal{D}$ if $\left(\mathcal{C}  - \mathcal{D} \right)$ is a TPSD (TPD) T-product tensor.

We have the following theorem from Theorem 5 in~\cite{zheng2021t}.
\begin{theorem}\label{thm:tensor and matrix pd relation}
If a T-product tensor $\mathcal{C} \in \mathbb{R}^{m \times m \times p}$ can be diagonalized as 
\begin{eqnarray}\label{eq:block diagonalized format}
\mbox{bcirc}(\mathcal{C}) = \left( \mathbf{F}^{\mathrm{H}}_m \otimes \mathbf{I}_m \right) \mbox{Diag}\left( \mathbf{C}_i: i \in \{1, \cdots, m \} \right)  \left( \mathbf{F}_m \otimes \mathbf{I}_m \right),
\end{eqnarray}
where $\mathbf{F}$ is the DFT matrix defined by Eq.~\eqref{eq:DFT def}; then $\mathcal{C}$ is symmetric, TPD (TPSD) if and only if all matrices $\mathbf{C}_i$ are Hermitian, positive definite (positive semidefinite).
\end{theorem}

Let $\mathcal{C} \in \mathbb{R}^{m \times m \times p}$ can be block diagonalized as Eq.~\eqref{eq:block diagonalized format}. Then, a real number $\lambda$ is said to be a \emph{T-eigenvalue} of $\mathcal{C}$, denoted as $\lambda ( \mathcal{C} )$, if it is an eigenvalue of some $\mathbf{C}_i$ for $i \in \{ 1, \cdots, m \}$. The largest and smallest T-eigenvalue of $\mathcal{C}$ are represented by $\lambda_{\max} (\mathcal{C} ) $ and $\lambda_{\min} (\mathcal{C} ) $, respectively. We use $\lambda_{i, j}$ for the $j$-th largest T-eigenvalue of the matrix $\mathbf{C}_i$.  We also use $\sigma_{i, j}$, named as \emph{T-singular values}, for the $j$-th largest singular values of the matrix $\mathbf{C}_i$. 

We define the T-product tensor \emph{trace} for a tensor $\mathcal{C} =  (c_{ijk}) \in \mathbb{C}^{m \times m \times p}$, denoted by $\mathrm{Tr}(\mathcal{C})$, as following
\begin{eqnarray}\label{eq:trace def}
\mathrm{Tr}(\mathcal{C}) \define \sum\limits_{i=1}^{m}\sum\limits_{k=1}^{p} c_{iik},
\end{eqnarray}
which is the summation of all entries in f-diagonal components. Then, we have the following lemma about trace properties. 
\begin{lemma}\label{lma:T product trace properties}
For any tensors $\mathcal{C}, \mathcal{D} \in \mathbb{C}^{m \times m \times p}$, we have 
\begin{eqnarray}\label{eq:trace linearity}
\mathrm{Tr}(c \mathcal{C} + d \mathcal{D}) = c  \mathrm{Tr}(\mathcal{C}) + d \mathrm{Tr}(\mathcal{D}),
\end{eqnarray}
where $c, d$ are two contants. And, the transpose operation will keep the same trace value, i.e., 
\begin{eqnarray}\label{eq:trace transpose same}
\mathrm{Tr}(\mathcal{C}) = \mathrm{Tr}(\mathcal{C}^{T}).
\end{eqnarray}
Finally, we have 
\begin{eqnarray}\label{eq:trace commutativity}
\mathrm{Tr}(\mathcal{C} \star  \mathcal{D} ) = \mathrm{Tr}(\mathcal{D} \star  \mathcal{C}).
\end{eqnarray}
\end{lemma}
\textbf{Proof:}
Eqs.~\eqref{eq:trace linearity} and~\eqref{eq:trace transpose same} are true from trace definition directly. 

From T-product definition, the $i$-th frontal slice matrix of $\mathcal{D} \star  \mathcal{C}$ is 
\begin{eqnarray}\label{eq:i slice DC}
\mathbf{D}^{(i)} \mathbf{C}^{(1)}  + \mathbf{D}^{(i-1)} \mathbf{C}^{(2)} + \cdots +
\mathbf{D}^{(1)} \mathbf{C}^{(i)} + \mathbf{D}^{(m)} \mathbf{C}^{(i+1) } + \cdots
+ \mathbf{D}^{ (i+1)} \mathbf{C}^{(m)}, 
\end{eqnarray}
similarly, the $i$-th frontal slice matrix of $\mathcal{C} \star  \mathcal{D}$ is 
\begin{eqnarray}\label{eq:i slice CD}
\mathbf{C}^{(i)} \mathbf{D}^{(1)}  + \mathbf{C}^{(i-1)} \mathbf{D}^{(2)} + \cdots +
\mathbf{C}^{(1)} \mathbf{D}^{(i)} + \mathbf{C}^{(m)} \mathbf{D}^{(i+1)} + \cdots
+ \mathbf{C}^{ (i+1)}  \mathbf{D}^{(m)}. 
\end{eqnarray}
Because the matrix trace of Eq.~\eqref{eq:i slice DC} and the matrix trace of Eq.~\eqref{eq:i slice CD} are same for each slice $i$ due to linearity and invariant under cyclic permutations of matrix trace, we have Eq.~\eqref{eq:trace commutativity} by summing over all frontal matrix slices. $\hfill \Box$


Below, we will define the \emph{determinant} of a T-product tensor $\mathcal{C} \in \mathbb{R}^{m \times m \times p}$, represented by $\det ( \mathcal{C} )$, as
\begin{eqnarray}\label{eq:def T prod tensor determinant}
\det ( \mathcal{C} ) &=& \prod\limits_{i=1, j=1}^{i=m, j=p}\lambda_{i, j}. 
\end{eqnarray}

%
%

We have the following theorem from Theorem 6 in~\cite{zheng2021t} about symmetric T-product tensor decomposition. 
\begin{theorem}\label{thm:T eigenvalue decomp}
Every symmetric T-product tensor $\mathcal{C} \in \mathbb{R}^{m \times m \times p}$ can be factored as 
\begin{eqnarray}
\mathcal{C} = \mathcal{U}^{T} \star\mathcal{D} \star \mathcal{U},
\end{eqnarray}
where $\mathcal{U}$ is an orthogonal tensor, i.e., $\mathcal{U}^{T} \star \mathcal{U} = \mathcal{I}_{m, m, p}$, and $\mathcal{D}$ is a F-diagonal tensor, i.e., each frontal slice of $\mathcal{D}$ is a diagonal matrix, such that diagonal entries of $\left( \mathbf{F}_m \otimes \mathbf{I}_m \right) \mbox{bcirc}\left( \mathcal{D} \right) \left( \mathbf{F}^{\mathrm{H}}_m \otimes \mathbf{I}_m \right)$ are T-eigenvalues of $\mathcal{C}$. If $\mathcal{C}$ is a TPD (TPSD) tensor, then all of its T-eigenvalues are positive (nonnegative). 
\end{theorem}

From Theorem~\ref{thm:T eigenvalue decomp} and Lemma~\ref{lma:T product trace properties}, we have the fact that 
\begin{eqnarray}\label{eq:trace is eigen sum}
\mathrm{Tr}(\mathcal{C}) = \sum\limits_i \lambda_i ( \mathcal{C} ).
\end{eqnarray}

%

If a symmetric T-product tensor $\mathcal{C} \in \mathbb{R}^{m \times m \times p}$ can be expressed as the format shown by Eq.~\eqref{eq:block diagonalized format}, the T-eigenvalues of  $\mathcal{C}$ with respect to the matrix $\mathbf{C}_i$ are denoted as $\lambda_{i, k_i}$, where $1 \leq k_i \leq m$, and we assume that 
$\lambda_{i, 1} \geq \lambda_{i, 2} \geq \cdots \geq \lambda_{i, m}$ (including multiplicities). Then, $\lambda_{i, k_i}$ is the $k_i$-th largest T-eigenvalue associated to the matrix $\mathbf{C}_i$. If we sort all T-eigenvalues of $\mathcal{C}$ from the largest one to the smallest one, we use $\tilde{k}$, a smallest integer between 1 to $m \times p$ (inclusive) associated with $p$ given positive integers $k_1, k_2, \cdots, k_p$ that satisfies
\begin{eqnarray}\label{eq:tilde k}
\lambda_{\tilde{k}} = \min\limits_{1 \leq i \leq m} \lambda_{i, k_i},~\mbox{and}~\lambda_{\tilde{k}}  \geq \lambda_{i, k_i + 1},
\end{eqnarray}
and we set $\tilde{i}$ from $\lambda_{\tilde{k}} $ as
\begin{eqnarray}\label{eq:tilde i}
\tilde{i}= \arg \min\limits_{i }\left\{  \lambda_{\tilde{k}}= \lambda_{i, k_i} \right\}.
\end{eqnarray}
Then, we will have the following Courant-Fischer theorem for T-product tensors. 


\begin{theorem}\label{thm:Courant-Fischer T-product}
Given a symmetric T-product tensor $\mathcal{C} \in \mathbb{R}^{m \times m \times p}$ and $p$ positive integers $k_1, k_2, \cdots, k_p$ with $1 \leq k_i \leq m$, then we have 
\begin{eqnarray}
\lambda_{\tilde{k}} &=& \max\limits_{\substack{S \in \mathbb{R}^{m \times 1 \times p}\\ \dim(\mathrm{S})  = \{k_1, \cdots, k_p \}   }} \min\limits_{\mathcal{X} \in S } \frac{ \langle \mathcal{X}, \mathcal{C} \star \mathcal{X} \rangle }{ \langle \mathcal{X}, \mathcal{X} \rangle }\nonumber \\ 
 &=&  \min\limits_{\substack{T \in \mathbb{R}^{m \times 1 \times p}\\ \dim(T)  = \{n- k_1, \cdots, n-k_{\tilde{i}-1},  n-k_{\tilde{i}}+1,  n-k_{\tilde{i}+1},\cdots,  n-k_p \}   }} \max\limits_{\mathcal{X} \in T} \frac{ \langle \mathcal{X}, \mathcal{C} \star \mathcal{X} \rangle }{ \langle \mathcal{X}, \mathcal{X} \rangle } 
\end{eqnarray}
where $\lambda_{\tilde{k}}$ and $\tilde{i}$ are defined by Eqs.~\eqref{eq:tilde k} and~\eqref{eq:tilde i}.
\end{theorem}
\textbf{Proof:}

First, we have to express $ \langle \mathcal{X}, \mathcal{C} \star \mathcal{X} \rangle $ by matrices of $\mathbf{C}_i$ and $\mathbf{X}_i$ through the representation shown by Eq.~\eqref{eq:block diagonalized format}. It is
\begin{eqnarray}\label{eq:thm:Courant-Fischer T-product}
\langle \mathcal{X}, \mathcal{C} \star \mathcal{X} \rangle &=& \frac{1}{p} \langle \mbox{bcirc}(\mathcal{X}), \mbox{bcirc}(\mathcal{C}) \mbox{bcirc}(\mathcal{X})    \rangle \nonumber \\
&=& \frac{1}{p} \mathrm{Tr} \left( \mbox{bcirc}(\mathcal{X})^{\mathrm{H}} \mbox{bcirc}(\mathcal{C}) \mbox{bcirc}(\mathcal{X}) \right)\nonumber \\
&=& \frac{1}{p} \mathrm{Tr} \left(\mathbf{F}^{\mathrm{H}}_p  \mbox{Diag}\left( \mathbf{x}^{\mathrm{H}}_i \mathbf{A}_i  \mathbf{x}_i: i \in \{1,\cdots,p\} \right)\mathbf{F}_p  \right)\nonumber \\
&=& \frac{1}{p} \mathrm{Tr} \left(  \mbox{Diag}\left( \mathbf{x}^{\mathrm{H}}_i \mathbf{A}_i  \mathbf{x}_i: i \in \{1,\cdots,p\} \right)  \right) = \frac{1}{p}\sum\limits_{i=1}^p \mathbf{x}^{\mathrm{H}}_i \mathbf{A}_i  \mathbf{x}_i
\end{eqnarray}

We will just verify the first characterization of $\lambda_{\tilde{k}}$. The other is similar. Let $S_i$ be the projection of $S$ to the space with dimension $k_i$ spanned by $\mathbf{v}_{i, 1}, \cdots, \mathbf{v}_{i, k_i}$, for every $\mathbf{x}_i \in S_{i}$, we can write $\mathbf{x}_i = \sum\limits^{k_i}_{j=1} c_{i, j}  \mathbf{v}_{i, j}$. To show that the value $\lambda_{\tilde{k}}$ is achievable, note that 
\begin{eqnarray}
 \frac{ \langle \mathcal{X}, \mathcal{C} \star \mathcal{X} \rangle }{ \langle \mathcal{X}, \mathcal{X} \rangle }
&=&  \frac{  \frac{1}{p}\sum\limits_{i=1}^p \mathbf{x}^{\mathrm{H}}_i \mathbf{A}_i  \mathbf{x}_i  }{  \frac{1}{p}\sum\limits_{i=1}^p \mathbf{x}^{\mathrm{H}}_i  \mathbf{x}_i  } = \frac{  \sum\limits_{i=1}^p    \sum\limits^{k_i}_{j=1} \lambda_{i, j}    c_{i, j}^{\ast} c_{i, j}          }{   \sum\limits_{i=1}^p    \sum\limits^{k_i}_{j=1}   c_{i, j}^{\ast} c_{i, j}  } \nonumber \\
&\geq & \frac{  \sum\limits_{i=1}^p    \sum\limits^{k_i}_{j=1} \lambda_{\tilde{k}}  c_{i, j}^{\ast} c_{i, j}          }{   \sum\limits_{i=1}^p    \sum\limits^{k_i}_{j=1}   c_{i, j}^{\ast} c_{i, j}  } = \lambda_{\tilde{k}}
\end{eqnarray}
To verify that this is the maximum, let $T_{\tilde{i}}$ be the projection of $T$ to the space with dimension $k_{\tilde{i}}$ with dimension $n- k_{\tilde{i}} + 1$, then the intersection of $S$ and $T_{\tilde{i}}$ is not empty. We have
\begin{eqnarray}
\min\limits_{\mathcal{X} \in S } \frac{ \langle \mathcal{X}, \mathcal{C} \star \mathcal{X} \rangle }{ \langle \mathcal{X}, \mathcal{X} \rangle } &\leq & \min\limits_{\mathcal{X} \in S \cap T} \frac{ \langle \mathcal{X}, \mathcal{C} \star \mathcal{X} \rangle }{ \langle \mathcal{X}, \mathcal{X} \rangle }.
\end{eqnarray}
Any such $\mathbf{x}_{\tilde{i}} \in S \cap T_{\tilde{i}}$ can be expressed as $\mathbf{x}_{\tilde{i}} =  \sum\limits^{m}_{ j=k_{\tilde{i} }} c_{\tilde{i}, j}  \mathbf{v}_{\tilde{i} j}$, and any $i$ for $i \neq \tilde{i}$, we have $\mathbf{x}_{i} \in S \cap T_{i}$ expressed as $\mathbf{x}_{ i } =  \sum\limits^{m}_{ j=k_{i} + 1} c_{i, j}  \mathbf{v}_{i, j}$. Then, we have 
\begin{eqnarray}
 \frac{ \langle \mathcal{X}, \mathcal{C} \star \mathcal{X} \rangle }{ \langle \mathcal{X}, \mathcal{X} \rangle }
&=&  \frac{  \frac{1}{p}\sum\limits_{i=1}^p \mathbf{x}^{\mathrm{H}}_i \mathbf{A}_i  \mathbf{x}_i  }{  \frac{1}{p}\sum\limits_{i=1}^p \mathbf{x}^{\mathrm{H}}_i  \mathbf{x}_i  } = \frac{  \sum\limits_{i=1}^p     \sum\limits^{m}_{ \substack{j=k_i + 1; i \neq \tilde{i} \\ j=k_{\tilde{i}};  i = \tilde{i}    }}   \lambda_{i, j}    c_{i, j}^{\ast} c_{i, j}          }{    \sum\limits_{i=1}^p    \sum\limits^{m}_{ \substack{j=k_i + 1; i \neq \tilde{i} \\ j=k_{\tilde{i}};  i = \tilde{i}    }   } c_{i, j}^{\ast} c_{i, j}   } \nonumber \\
&\leq & \frac{  \sum\limits_{i=1}^p    \sum\limits^{m}_{ \substack{j=k_i + 1; i \neq \tilde{i} \\ j=k_{\tilde{i}};  i = \tilde{i}    }   } \lambda_{\tilde{k}}  c_{i, j}^{\ast} c_{i, j}          }{   \sum\limits_{i=1}^p    \sum\limits^{m}_{ \substack{j=k_i + 1; i \neq \tilde{i} \\ j=k_{\tilde{i}};  i = \tilde{i}    }   } c_{i, j}^{\ast} c_{i, j}  } = \lambda_{\tilde{k}}.
\end{eqnarray}
Therefore, for all subspaces $S$ of dimensions $\{k_1, \cdots, k_p\}$, we have $\min\limits_{\mathcal{X} \in S} \frac{ \langle \mathcal{X}, \mathcal{C} \star \mathcal{X} \rangle }{ \langle \mathcal{X}, \mathcal{X} \rangle } \leq \lambda_{\tilde{k}}$
$\hfill \Box$

Given a symmetric T-product tensor $\mathcal{C}$ with associated matrices $\mathbf{C}_i$ provided by Eq.~\eqref{eq:block diagonalized format}, next theorem is the representation of the summation of all the largest $k_i$ T-eigenvalues of $\mathbf{C}_i$ and the summation of all the smallest $k_i$ T-eigenvalues of $\mathbf{C}_i$. 
\begin{theorem}\label{thm:extreme T-eigenvalues sum rep}
Let $\mathcal{C} \in \mathbb{R}^{m \times m \times p}$ be a symmetric T-product tensor with associated matrices $\mathbf{C}_i$ provided by Eq.~\eqref{eq:block diagonalized format}, and we sort T-eigenvalues of the matrix $\mathbf{C}_i$ as $\lambda_{i, 1} \geq \lambda_{i, 2} \geq \cdots \geq \lambda_{i, k_i}$. Then, we have
\begin{eqnarray}\label{eq:largest k i T eigenvalues}
\sum\limits_{i=1}^{p} \max\limits_{\mathbf{U}_i \mathbf{U}^{\mathrm{H}}_i = \mathbf{I}_{k_i} } \mathrm{Tr} 
\left( \mathbf{U}_i \mathbf{C}_i \mathbf{U}^{\mathrm{H}}_i \right) &=& 
\sum\limits_{i=1}^{p} \sum\limits_{j = 1}^{k_i} \lambda_{i, j}(\mathbf{C}_i);
\end{eqnarray}
and
\begin{eqnarray}\label{eq:smallest k i T eigenvalues}
\sum\limits_{i=1}^{p} \min\limits_{\mathbf{U}_i \mathbf{U}^{\mathrm{H}}_i = \mathbf{I}_{k_i} } \mathrm{Tr} 
\left(  \mathbf{U}_i \mathbf{C}_i \mathbf{U}^{\mathrm{H}}_i \right) &=& 
\sum\limits_{i=1}^{p} \sum\limits_{j = 1}^{k_i} \lambda_{i, m-j+1}(\mathbf{C}_i),
\end{eqnarray}
where $\mathbf{U}_i$ are $k_{i} \times m$ complex matrices. 
\end{theorem}
\textbf{Proof:}
From Theorem~\ref{thm:tensor and matrix pd relation}, we may assume that $\mathbf{C}_i$ are diagonal matrices, denoted as $\mathbf{D}_i$, since $\mathbf{C}_i$ are symmetric T-product matrices. Therefore, we have the expression $\mathbf{C}_i = \mathbf{V} \mathbf{D}_i \mathbf{V}^{\mathrm{H}}$. Then, we have
\begin{eqnarray}\label{eq1:thm:extreme T-eigenvalues sum rep}
\mathrm{Tr} \left( \mathbf{U}_i \mathbf{D}_i \mathbf{U}^{\mathrm{H}}_i  \right) = 
\sum\limits_{j=1}^{k_i}  \sum\limits_{l=1}^{m}   u_{j, l}^{\ast} u_{j, l} \lambda_{i, l}(\mathbf{C}_i)
= 
\sum\limits_{j=1}^{k_i}  \sum\limits_{l=1}^{m}   p_{j, l} \lambda_{i, l}(\mathbf{C}_i)
= [\overbrace{1,1,\cdots,1}^{\mbox{$k_i$ terms}}] \mathbf{P}  \left[
    \begin{array}{c}
       \lambda_{i,1}  \\
       \lambda_{i,2}  \\
       \vdots    \\
       \lambda_{i,m}  \\
    \end{array}
\right],
\end{eqnarray} 
where $\mathbf{P} = (p_{j,l})$ is a $k_i \times m$ stochastic matrix. Then, we can concatenate an $(m - k_i) \times m$ matrix $\mathbf{Q}$ to the matrix $\mathbf{P}$ to make the following matrix $\left[
    \begin{array}{c}
\mathbf{P}  \\
\mathbf{Q}  \\
    \end{array}
\right]$ as doubly  stochastic from 2.C.1(4) from~\cite{marshall2011inequalities}. Then, Eq.~\eqref{eq1:thm:extreme T-eigenvalues sum rep} can be  expressed as
\begin{eqnarray}\label{eq2:thm:extreme T-eigenvalues sum rep}
\mathrm{Tr} \left( \mathbf{U}_i \mathbf{D}_i \mathbf{U}^{\mathrm{H}}_i  \right) =  
[\overbrace{1,1,\cdots,1}^{\mbox{$k_i$ terms}},\overbrace{0,0,\cdots,0}^{\mbox{$m - k_i$ terms}}] \left[
    \begin{array}{c}
\mathbf{P}  \\
\mathbf{Q}  \\
    \end{array}
\right]  \left[
    \begin{array}{c}
       \lambda_{i,1}  \\
       \lambda_{i,2}  \\
       \vdots    \\
       \lambda_{i,m}  \\
    \end{array}
\right]
\end{eqnarray} 

Given two lists of real numbers, $[a_1, \cdots, a_n]$ and $[b_1, \cdots, b_n]$, we use $[a_1, \cdots, a_n] \prec [b_1, \cdots, b_n]$ to represent the following relationships:
\begin{eqnarray}
\sum\limits_{i=1}^{k} a_i \leq \sum\limits_{i=1}^{k} b_i, 
\end{eqnarray}
holds for any $k$ between $1$ and $n$. From Eq.~\eqref{eq2:thm:extreme T-eigenvalues sum rep}, we have $[\lambda_{i,1}, \cdots, \lambda_{i,m}] [\mathbf{P}^{\mathrm{T}},\mathbf{Q}^{\mathrm{T}}] \prec [\lambda_{i,1}, \cdots, \lambda_{i,m}] $ and 3.H.2.b from~\cite{marshall2011inequalities}, we have 
\begin{eqnarray}\label{eq3:thm:extreme T-eigenvalues sum rep}
\mathrm{Tr} 
\left( \mathbf{U}_i \mathbf{C}_i \mathbf{U}^{\mathrm{H}}_i \right) &\leq & 
 \sum\limits_{j = 1}^{k_i} \lambda_{i, j}(\mathbf{C}_i);
\end{eqnarray}
and
\begin{eqnarray}\label{eq4:thm:extreme T-eigenvalues sum rep}
\mathrm{Tr} 
\left(  \mathbf{U}_i \mathbf{C}_i \mathbf{U}^{\mathrm{H}}_i \right) &\geq & 
 \sum\limits_{j = 1}^{k_i} \lambda_{i, m-j+1}(\mathbf{C}_i).
\end{eqnarray}
Finally, this theorem is proved by applying $\sum\limits_{i=1}^{p}$ to both sides of Eqs.~\eqref{eq3:thm:extreme T-eigenvalues sum rep} and~\eqref{eq4:thm:extreme T-eigenvalues sum rep} with respect to the index $i$, and note that $\mathbf{U}_i \mathbf{V}_i = \left( \mathbf{I}_{k_i}, \mathbf{O} \right)$ and $\mathbf{U}_i \mathbf{V}_i = \left( \mathbf{O}, \mathbf{I}_{k_i} \right)$, respectively.
$\hfill \Box$


\subsection{Unitarily Invariant T-product Tensor Norms}\label{sec:Unitarily Invariant T-product Tensor Norms}


Let us represent the T-eigenvalues of a symmetric T-product tensor $\mathcal{H} \in \mathbb{R}^{m \times m \times p} $ in decreasing order by the vector $\vec{\lambda}(\mathcal{H}) = (\lambda_1(\mathcal{H}), \cdots, \lambda_{m \times p}(\mathcal{H}))$, where $m \times p$ is the total number of T-eigenvalues. We use $\mathbb{R}_{\geq 0} (\mathbb{R}_{> 0})$ to represent a set of nonnegative (positive) real numbers. Let $\left\Vert \cdot \right\Vert_{\rho}$ be a unitarily invariant tensor norm, i.e., $\left\Vert \mathcal{H}\star \mathcal{U}\right\Vert_{\rho} = \left\Vert \mathcal{U}\star \mathcal{H}\right\Vert_{\rho} = \left\Vert \mathcal{H}\right\Vert_{\rho} $,  where $\mathcal{U}$ is any unitary tensor. Let $\rho : \mathbb{R}_{\geq 0}^{m \times p} \rightarrow \mathbb{R}_{\geq 0}$ be the corresponding gauge function that satisfies H$\ddot{o}$lder’s inequality so that 
\begin{eqnarray}\label{eq:def gauge func and general unitarily invariant norm}
\left\Vert \mathcal{H} \right\Vert_{\rho} = \left\Vert |\mathcal{H}| \right\Vert_{\rho} = \rho(\vec{\lambda}( | \mathcal{H} | ) ),
\end{eqnarray}
where $ |\mathcal{H}|  \define \sqrt{\mathcal{H}^H \star \mathcal{H}} $. The bijective correspondence between symmetric gauge functions on $\mathbb{R}_{\geq 0}^{m \times p}$ and unitarily invariant norms is due to von Neumann~\cite{fan1955some}. 

Several popular norms can be treated as special cases of unitarily invariant tensor norm. The first one is Ky Fan like $k$-norm~\cite{fan1955some} for tensors. For $k \in \{1,2,\cdots,m \times p\}$, the Ky Fan $k$-norm~\cite{fan1955some} for tensors  $\mathcal{H}  \mathbb{R}^{m \times m \times p} $, denoted as $\left\Vert \mathcal{H}\right\Vert_{(k)}$, is defined as:
\begin{eqnarray}\label{eq: Ky Fan k norm for tensors}
\left\Vert \mathcal{H}\right\Vert_{(k)} \define \sum\limits_{i=1}^{k} \lambda_i(  |\mathcal{H}|  ).
\end{eqnarray}
If $k=1$,  the Ky Fan $k$-norm for tensors is the tensor operator norm, denoted as $ \left\Vert \mathcal{H} \right\Vert$. The second one is Schatten $p$-norm for tensors, denoted as $\left\Vert \mathcal{H}\right\Vert_{p}$, is defined as:
\begin{eqnarray}\label{eq: Schatten p norm for tensors}
\left\Vert \mathcal{H}\right\Vert_{p} \define (\mathrm{Tr}|\mathcal{H}|^p )^{\frac{1}{p}},
\end{eqnarray}
where $ p \geq 1$. If $p=1$, it is the trace norm. 

Following inequality is the extension of H\"{o}lder inequality to gauge function $\rho$ which will be used later to prove majorization relations. 
\begin{lemma}\label{lma:Holder inquality for gauge function}
For $n$ nonnegative real vectors with the dimension $r$, i.e., $\mathbf{b}_i = (b_{i_1}, \cdots, b_{i_r}) \in \mathbb{R}_{\geq 0}^{r}$, and $\alpha > 0$ with $\sum\limits_{i=1}^n \alpha_i = 1$, we have 
\begin{eqnarray}\label{eq1:lma:Holder inquality for gauge function}
\rho\left( \prod\limits_{i=1}^n b_{i_1}^{\alpha_i},  \prod\limits_{i=1}^n b_{i_2}^{\alpha_i}, \cdots,  \prod\limits_{i=1}^n b_{i_r}^{\alpha_i}  \right) \leq  \prod\limits_{i=1}^n \rho(\mathbf{b}_i)^{\alpha_i} 
\end{eqnarray}
\end{lemma}
\textbf{Proof:}
This proof is based on mathematical induction. The base case for $n=2$ has been shown by Theorem IV.1.6 from~\cite{bhatia2013matrix}. 

We assume that Eq.~\eqref{eq1:lma:Holder inquality for gauge function} is true for $n=m$, where $m > 2$. Let $\odot$ be the component-wise product (Hadamard product) between two vectors.  Then, we have 
\begin{eqnarray}\label{eq2:lma:Holder inquality for gauge function}
\rho\left( \prod\limits_{i=1}^{m+1} b_{i_1}^{\alpha_i},  \prod\limits_{i=1}^{m+1} b_{i_2}^{\alpha_i}, \cdots,  \prod\limits_{i=1}^{m+1} b_{i_r}^{\alpha_i}  \right) = 
\rho\left( \odot_{i=1}^{m+1} \mathbf{b}_i^{\alpha_i}  \right),
\end{eqnarray}
where $\odot_{i=1}^{m+1} \mathbf{b}_i^{\alpha_i}$ is defined as $\left( \prod\limits_{i=1}^{m+1} b_{i_1}^{\alpha_i},  \prod\limits_{i=1}^{m+1} b_{i_2}^{\alpha_i}, \cdots,  \prod\limits_{i=1}^{m+1} b_{i_r}^{\alpha_i}  \right)$ with $\mathbf{b}_i^{\alpha_i} \define (b_{i_1}^{\alpha_i}, \cdots, b_{i_r}^{\alpha_i})$. Under such notations, Eq.~\eqref{eq2:lma:Holder inquality for gauge function} can be bounded as  
\begin{eqnarray}\label{eq3:lma:Holder inquality for gauge function}
\rho\left( \odot_{i=1}^{m+1} \mathbf{b}_i^{\alpha_i}  \right) &= &
\rho\left( \left( \odot_{i=1}^{m} \mathbf{b}_i^{\frac{\alpha_i}{ \sum\limits_{j=1}^m \alpha_j }  } \right)^{\sum\limits_{j=1}^m \alpha_j } \odot \mathbf{b}_{m+1}^{\alpha_{m+1}}\right) \nonumber \\
& \leq & \left[ \rho^{\sum\limits_{j=1}^m \alpha_j } \left( \odot_{i=1}^{m} \mathbf{b}_i^{\frac{\alpha_i}{ \sum\limits_{j=1}^m \alpha_j }  }  \right)  \right] \cdot \rho( \mathbf{b}_{m+1})^{\alpha_{m+1}} \leq  \prod\limits_{i=1}^{m+1} \rho(\mathbf{b}_i)^{\alpha_i}. 
\end{eqnarray}
By mathematical induction, this lemma is proved. $\hfill \Box$

\subsection{Antisymmetric Kronecker Product for T-product Tensors}\label{sec:Antisymmetric  Kronecker Product for T-product Tensors}

In this section, we will discuss a machinery of antisymmetric Kronecker product for T-product tensors and this scheme will be used later for log-majorization results. Let $\mathfrak{H}$ be an $m \times p$-dimensional Hilbert space. For each $k \in \mathbb{N}$, let $\mathfrak{H}^{\otimes k}$ denote the $k$-fold Kronecker product of $\mathfrak{H}$, which is the $(m \times p)^k$-dimensional Hilbert space with respect to the inner product defined by
\begin{eqnarray}
\langle \mathbf{X}_1 \otimes \cdots \otimes \mathbf{X}_k,  \mathbf{Y}_1 \otimes \cdots \otimes \mathbf{Y}_k \rangle \define \prod\limits_{i=1}^k \langle \mathbf{X}_i, \mathbf{Y}_i \rangle.
\end{eqnarray}
For $\mathbf{X}_1, \cdots, \mathbf{X}_k \in \mathfrak{H}$, we define $\mathbf{X}_1 \wedge \cdots \wedge  \mathbf{X}_k \in \mathfrak{H}^{\otimes k}$ by
\begin{eqnarray}
\mathbf{X}_1 \wedge \cdots \wedge  \mathbf{X}_k \define \frac{1}{\sqrt{k!}} \sum\limits_{\sigma} (\mbox{sgn} \sigma) \mathbf{X}_{\sigma(1)} \otimes \cdots \otimes  \mathbf{X}_{\sigma(k)}, 
\end{eqnarray}
where $\sigma$ runs over all permutations on $\{1, 2, \cdots, k\}$ and $\mbox{sgn} \sigma = \pm 1$ depending on $\sigma$ is even or odd. The subspace of $\mathfrak{H}^{\otimes k}$ spanned by $\{\mathbf{X}_1 \wedge \cdots \wedge  \mathbf{X}_k \}$, where $\mathbf{X}_i \in \mathfrak{H}$, is named as $k$-fold antisymmetric Kronecker product of $\mathfrak{H}$ and represented by $\mathfrak{H}^{\wedge k}$. 

For each $\mathcal{C} \in \mathbb{R}^{m \times m \times p}$ and $k \in \mathbb{N}$, the $k$-fold Kronecker product $\mathcal{C}^{\otimes k} \in \mathbb{R}^{m^k \times m^k \times p^k} $ is given by
\begin{eqnarray}
\mathcal{C}^{\otimes k} \star \left( \mathbf{X}_1 \otimes \cdots \otimes \mathbf{X}_k \right) \define 
\left( \mathcal{C} \star \mathbf{X}_1 \right) \otimes  \cdots \otimes \left( \mathcal{C} \star \mathbf{X}_k \right).
\end{eqnarray}
Because $\mathfrak{H}^{\wedge k}$ is invariant for $\mathcal{C}^{\otimes k}$, the antisymmetric Kronecker product of $\mathcal{C}^{\wedge k}$ of $\mathcal{C}$ can be defined as $\mathcal{C}^{\wedge k} = \mathcal{C}^{\otimes}|_{\mathfrak{H}^{\wedge k}}$, then we have
\begin{eqnarray}
\mathcal{C}^{\wedge k} \star \left( \mathbf{X}_1 \wedge \cdots \wedge \mathbf{X}_k \right) 
= \left( \mathcal{C} \star \mathbf{X}_1 \right)  \wedge  \cdots \wedge  \left( \mathcal{C} \star \mathbf{X}_k \right).
\end{eqnarray}

We will provide the following lemmas about antisymmetric Kronecker product.
\begin{lemma}\label{lma:antisymmetric tensor product properties}
Let $\mathcal{A}, \mathcal{B}, \mathcal{C},  \mathcal{E} \in \mathbb{R}^{m \times m \times p}$ be T-product tensors , for any $k \in \{1,2,\cdots,m \times p\}$, we have 
\begin{enumerate}
	\item $(\mathcal{A}^{\wedge k})^\mathrm{T} = (\mathcal{A}^\mathrm{T} )^{\wedge k}$.
	\item $(\mathcal{A}^{\wedge k}) \star (\mathcal{B}^{\wedge k})= (\mathcal{A}\star \mathcal{B})^{\wedge k}$. 
	\item If $\lim\limits_{i \rightarrow \infty} \left\Vert \mathcal{A}_i -  \mathcal{A} \right\Vert \rightarrow 0$ , then $\lim\limits_{i \rightarrow \infty} \left\Vert \mathcal{A}^{\wedge k}_i -  \mathcal{A}^{\wedge k} \right\Vert \rightarrow 0$.
	\item If $\mathcal{C} \succeq \mathcal{O}$ (zero tensor), then $\mathcal{C}^{\wedge k} \succeq \mathcal{O}$ and $(\mathcal{C}^p)^{\wedge k} = (\mathcal{C}^{\wedge k})^p$ for all $p \in \mathbb{R}_{> 0 }$.
    \item  $|\mathcal{A}|^{\wedge k} = | \mathcal{A}^{\wedge k}|$.
    \item  If $\mathcal{E} \succeq \mathcal{O}$ and $\mathcal{E}$ is invertible,  $(\mathcal{E}^z)^{\wedge k} = (\mathcal{E}^{\wedge k})^z$ for all $z \in \mathbb{E}$.
    \item  $\left\Vert \mathcal{E}^{\wedge k} \right\Vert = \prod\limits_{i=1}^{k} \lambda_i ( | \mathcal{E} |)$.
\end{enumerate}
\end{lemma}
\textbf{Proof:}
Items 1 and 2 are the restrictions of the associated relations $(\mathcal{A}^H)^{\otimes k} = (\mathcal{A}^{\otimes k})^H$ and $(\mathcal{A} \star \mathcal{B})^{\otimes k} = (\mathcal{A}^{\otimes k})\star (\mathcal{B}^{\otimes k})$ to $\mathfrak{H}^{\wedge k}$. The item 3 is true since, if $\lim\limits_{i \rightarrow \infty} \left\Vert \mathcal{A}_i -  \mathcal{A} \right\Vert \rightarrow 0$, we have $\lim\limits_{i \rightarrow \infty} \left\Vert \mathcal{A}^{\otimes k}_i -  \mathcal{A}^{\otimes k} \right\Vert \rightarrow 0$ and the asscoaited restrictions of $\mathcal{A}_i^{\otimes k}, \mathcal{A}^{\otimes k}$ to the antisymmetric subspace $\mathfrak{H}^k$. 

For the item 4, if $\mathcal{C} \succeq \mathcal{O}$, then we have $\mathcal{C}^{\wedge k} = ((\mathcal{C}^{1/2})^{\wedge k})^H \star  ((\mathcal{C}^{1/2})^{\wedge k}) \succeq   \mathcal{O}$ from items 1 and 2. If $p$ is ratonal, we have  $(\mathcal{C}^p)^{\wedge k} = (\mathcal{C}^{\wedge k})^p$  from the item 2, and the equality $(\mathcal{C}^p)^{\wedge k} = (\mathcal{C}^{\wedge k})^p$ is also true for any $p > 0$ if we apply the item 3 to approximate any irrelational numbers by rational numbers.

Because we have 
\begin{eqnarray}
|\mathcal{A}|^{\wedge k} =  \left( \sqrt{\mathcal{A}^H \mathcal{A}}\right)^{\wedge k}  =   \sqrt{ (\mathcal{A}^{\wedge k})^H \mathcal{A}^{\wedge k}  }=  | \mathcal{A}^{\wedge k}|,
\end{eqnarray}
from items 1, 2 and 4, so the item 5 is valid. 

For item 6, if $z  < 0$, item 6 is true for all $z \in \mathbb{R}$ by applying the item 4 to $\mathcal{E}^{-1}$. Since we can apply the definition $\mathcal{E}^{z} \define \exp(z \ln \mathcal{E})$ to have
\begin{eqnarray}
\mathcal{C}^p &=& \mathcal{E}^z~~\leftrightarrow~~\mathcal{C} = \exp\left(\frac{z}{p} \ln \mathcal{E} \right),
\end{eqnarray}
where $\mathcal{C} \succeq \mathcal{O}$. The general case of any $z \in \mathbb{C}$ is also true by applying the item 4 to $\mathcal{C} = \exp(\frac{z}{p} \ln \mathcal{E})$. 

For the item 7 proof, it is enough to prove the case that $\mathcal{E} \succeq \mathcal{O}$ due to the item 5. Then, from Theorem~\ref{thm:T eigenvalue decomp}, there exists a set of orthogonal tensors $\{\mathcal{U}_1, \cdots, \mathcal{U}_r\}$ such that $ | \mathcal{E} | \star \mathcal{U}_i = \lambda_i   \mathcal{U}_i$ for $1 \leq i \leq m \times p$. We then have 
\begin{eqnarray}
 | \mathcal{E} |^{\wedge k} \left(\mathcal{U}_{i_1} \wedge \cdots \wedge \mathcal{U}_{i_k}\right)
&=&  | \mathcal{E} |\star  \mathcal{U}_{i_1} \wedge \cdots \wedge  | \mathcal{E} |\star  \mathcal{U}_{i_k}  \nonumber \\
&=& \left( \prod\limits_{i=1}^{k} \lambda_i (  | \mathcal{E} |)  \right) \mathcal{U}_{i_1} \wedge \cdots \wedge \mathcal{U}_{i_k},
\end{eqnarray}
where $1 \leq i_1 < i_2 < \cdots < i_k \leq m \times p$. Hence, $\left\Vert  | \mathcal{E} |^{\wedge k} \right\Vert = \prod\limits_{i=1}^{k} \lambda_i (  | \mathcal{E} |)$. 
$\hfill \Box$

\section{Multivariate T-product Tensor Norm Inequalities}\label{sec:Multivariate T-product Tensor Norm Inequalities}

In this section, we will begin with the introduction of majorization techniques in Section~\ref{sec:Majorization Basis}.
Then, the majorization with integral average and log-majorization with integral average will be introduced by Section~\ref{sec:Majorization wtih Integral Average} and Section~\ref{sec:Log-Majorization wtih Integral Average}. These majorization results will be used to prove T-product tensor norm inequalities in Section~\ref{sec:T-product Tensor Norm Inequalities by Majorization}.


\subsection{Majorization Basis}\label{sec:Majorization Basis} 

In this subsection, we will discuss majorization and several lemmas about majorization which will be used at later proofs. 

Let $\mathbf{x} = [x_1, \cdots,x_r] \in \mathbb{R}^{m \times p}, \mathbf{y} = [y_1, \cdots,y_r] \in \mathbb{R}^{m \times p}$ be two vectors with following orders among entries $x_1 \geq \cdots \geq x_r$ and $y_1 \geq \cdots \geq y_r$, \emph{weak majorization} between vectors $\mathbf{x}, \mathbf{y}$, represented by $\mathbf{x} \prec_{w} \mathbf{y}$, requires following relation for  vectors $\mathbf{x}, \mathbf{y}$:
\begin{eqnarray}\label{eq:weak majorization def}
\sum\limits_{i=1}^k x_i \leq \sum\limits_{i=1}^k y_i,
\end{eqnarray}
where $k \in \{1,2,\cdots,r\}$. \emph{Majorization} between vectors $\mathbf{x}, \mathbf{y}$, indicated by $\mathbf{x} \prec \mathbf{y}$, requires following relation for vectors $\mathbf{x}, \mathbf{y}$:
\begin{eqnarray}\label{eq:majorization def}
\sum\limits_{i=1}^k x_i &\leq& \sum\limits_{i=1}^k y_i,~~\mbox{for $1 \leq k < r$;} \nonumber \\
\sum\limits_{i=1}^{m \times p} x_i &=& \sum\limits_{i=1}^{m \times p} y_i,~~\mbox{for $k = r$.}
\end{eqnarray}

For $\mathbf{x}, \mathbf{y} \in \mathbb{R}^{m \times p}_{\geq 0}$ such that  $x_1 \geq \cdots \geq x_r$ and $y_1 \geq \cdots \geq y_r$,  \emph{weak log majorization} between vectors $\mathbf{x}, \mathbf{y}$, represented by $\mathbf{x} \prec_{w \log} \mathbf{y}$, requires following relation for vectors $\mathbf{x}, \mathbf{y}$:
\begin{eqnarray}\label{eq:weak log majorization def}
\prod\limits_{i=1}^k x_i \leq \prod\limits_{i=1}^k y_i,
\end{eqnarray}
where $k \in \{1,2,\cdots,r\}$, and \emph{log majorization} between vectors $\mathbf{x}, \mathbf{y}$, represented by $\mathbf{x} \prec_{\log} \mathbf{y}$, requires equality for $k=r$ in Eq.~\eqref{eq:weak log majorization def}. If $f$ is a single variable function, $f(\mathbf{x})$ represents a vector of $[f(x_1),\cdots,f(x_r)]$. From Lemma 1 in~\cite{hiai2017generalized}, we have 
\begin{lemma}\label{lma:Lemma 1 Gen Log Hiai}
(1) For any convex function $f: [0, \infty) \rightarrow [0, \infty)$, if we have $\mathbf{x} \prec \mathbf{y}$, then $f(\mathbf{x}) \prec_{w} f(\mathbf{y})$. \\
(2) For any convex function and non-decreasing $f: [0, \infty) \rightarrow [0, \infty)$, if we have $\mathbf{x} \prec_{w} \mathbf{y}$, then $f(\mathbf{x}) \prec_{w} f(\mathbf{y})$. \\
\end{lemma}

Another lemma is from Lemma 12 in~\cite{hiai2017generalized}, we have 
\begin{lemma}\label{lma:Lemma 12 Gen Log Hiai}
Let $\mathbf{x}, \mathbf{y} \in \mathbb{R}^{m \times p}_{\geq 0}$ such that  $x_1 \geq \cdots \geq x_r$ and $y_1 \geq \cdots \geq y_r$ with $\mathbf{x}\prec_{\log} \mathbf{y}$. Also let $\mathbf{y}_i = [y_{i;1}, \cdots , y_{i;r} ] \in \mathbb{R}^{m \times p}_{\geq 0}$ be a sequence of vectors such that $y_{i;1} \geq \cdots \geq y_{i;r} > 0$ and $\mathbf{y}_i \rightarrow \mathbf{y}$ as $i \rightarrow \infty$. Then, there exists $i_0 \in \mathbb{N}$ and $\mathbf{x}_i  = [x_{i;1}, \cdots , x_{i;r} ] \in \mathbb{R}^{m \times p}_{\geq 0}$ for $i \geq i_0$ such that $x_{i;1} \geq \cdots \geq x_{i;r} > 0$, $\mathbf{x}_i \rightarrow \mathbf{x}$ as $i \rightarrow \infty$, and 
\begin{eqnarray}
\mathbf{x}_i \prec_{\log} \mathbf{y}_i \mbox{~~for $i \geq i_0$.} 
\end{eqnarray}
\end{lemma}

For any function $f$ on $\mathbb{R}_{\geq 0}$, the term $f(\mathbf{x}$ is defined as $f(\mathbf{x}) \define (f(x_1), \cdots, f(x_r))$ with conventions $e^{ - \infty} = 0$ and $\log 0 = - \infty$. 

\subsection{Majorization with Integral Average}\label{sec:Majorization wtih Integral Average}

Let $\Omega$ be a $\sigma$-compact metric space and $\nu$ a probability measure on the Borel $\sigma$-field of $\Omega$. Let $\mathcal{C}, \mathcal{D}_\tau \in \mathbb{R}^{m \times m \times p}$ be symmetric T-product tensors. We further assume that tensors $\mathcal{C}, \mathcal{D}_\tau$ are uniformly bounded in their norm for $\tau \in \Omega$. Let $\tau \in\Omega \rightarrow  \mathcal{D}_\tau$ be a continuous function such that $\sup \{\left\Vert  D_{\tau} \right\Vert: \tau \in \Omega  \} < \infty$. For notational convenience, we define the following relation:
\begin{eqnarray}\label{eq:integral eigen vector rep}
\left[ \int_{\Omega} \lambda_1(\mathcal{D}_\tau) d\nu(\tau), \cdots, \int_{\Omega} \lambda_{m \times p}(\mathcal{D}_\tau) d\nu(\tau) \right] \define \int_{\Omega^{m \times p}} \vec{\lambda}(\mathcal{D}_\tau) d\nu^{m \times p}(\tau).
\end{eqnarray}
If $f$ is a single variable function, the notation $f(\mathcal{C})$ represents a tensor function with respect to the tensor $\mathcal{C}$. 

\begin{theorem}\label{thm:weak int average thm 4}
Let $\Omega, \nu, \mathcal{C}, \mathcal{D}_\tau$ be defined as the beginning part of Section~\ref{sec:Majorization wtih Integral Average}, and $f: \mathbb{R} \rightarrow [0, \infty)$ be a non-decreasing convex function, we have following two equivalent statements:
\begin{eqnarray}\label{eq1:thm:weak int average thm 4}
\vec{\lambda}(\mathcal{C}) \prec_w  \int_{\Omega^{m \times p}} \vec{\lambda}(\mathcal{D}_\tau) d\nu^{m \times p}(\tau) \Longleftrightarrow \left\Vert f(\mathcal{C}) \right\Vert_{\rho} \leq 
\int_{\Omega} \left\Vert f(\mathcal{D}_{\tau}) \right\Vert_{\rho}  d\nu(\tau),
\end{eqnarray}
where $\left\Vert \cdot \right\Vert_{\rho}$ is the unitarily invariant norm defined in Eq.~\eqref{eq:def gauge func and general unitarily invariant norm}. 
\end{theorem}
\textbf{Proof:}
We assume that the left statement of Eq.~\eqref{eq1:thm:weak int average thm 4} is true and the function $f$ is a non-decreasing convex function. From Lemma~\ref{lma:Lemma 1 Gen Log Hiai}, we have 
\begin{eqnarray}\label{eq2:thm:weak int average thm 4}
\vec{\lambda}(f (\mathcal{C})) = f (\vec{\lambda}(\mathcal{C})) \prec_w  f \left(\int_{\Omega^{m \times p}} \vec{\lambda}(\mathcal{D}_\tau) d\nu^{m \times p}(\tau) \right).
\end{eqnarray}
From the convexity of $f$, we also have 
\begin{eqnarray}\label{eq3:thm:weak int average thm 4}
f \left(\int_{\Omega^{m \times p} } \vec{\lambda}(\mathcal{D}_\tau) d\nu^{m \times p}  (\tau) \right) \leq \int_{\Omega^{m \times p} } f(\vec{\lambda}(\mathcal{D}_\tau)) d\nu^{m \times p} (\tau) = \int_{\Omega^{m \times p}} \vec{\lambda} ( f(\mathcal{D}_\tau)) d\nu^{m \times p}(\tau).
\end{eqnarray}
Then, we obtain $\vec{\lambda}(f (\mathcal{C}))  \prec_{w} = \int_{\Omega^{m \times p}} \vec{\lambda} ( f(\mathcal{D}_\tau)) d\nu^{m \times p} (\tau)$. By applying Lemma 4.4.2 in~\cite{hiai2010matrix} to both sides of $\vec{\lambda}(f (\mathcal{C}))  \prec_{w} = \int_{\Omega^{m \times p} } \vec{\lambda} ( f(\mathcal{D}_\tau)) d\nu^{m \times p} (\tau)$ with gauge function $\rho$, we obtain 
\begin{eqnarray}\label{eq4:thm:weak int average thm 4}
\left \Vert f(\mathcal{C}) \right\Vert_{\rho} &\leq &\rho \left( \int_{\Omega^{m \times p} } \vec{\lambda} ( f(\mathcal{D}_\tau)) d\nu^{m \times p}  (\tau)  \right)  \nonumber \\
&\leq & \int_{\Omega} \rho(\vec{\lambda} ( f(\mathcal{D}_\tau))) d\nu(\tau) 
= \int_{\Omega} \left\Vert f(\mathcal{D}_\tau) \right\Vert_{\rho} d\nu(\tau).
\end{eqnarray}
Therefore, the right statement of Eq.~\eqref{eq1:thm:weak int average thm 4} is true from the left statement. 

On the other hand, if the right statement of Eq.~\eqref{eq1:thm:weak int average thm 4} is true, we select a function $f \define \max\{x + c, 0\} $, where $c$ is a positive real constant satisfying $\mathcal{C} + c \mathcal{I} \geq \mathcal{O}$, $\mathcal{D}_{\tau} + c \mathcal{I} \geq \mathcal{O}$ for all $\tau \in \Omega$, and tensors $\mathcal{C} + c \mathcal{I}, \mathcal{D}_{\tau} + c \mathcal{I}$. If the Ky Fan $k$-norm at the right statement of Eq.~\eqref{eq1:thm:weak int average thm 4} is applied, we have 
\begin{eqnarray}\label{eq5:thm:weak int average thm 4}
\sum\limits_{i=1}^k (\lambda_i (\mathcal{C}) + c ) \leq  
\sum\limits_{i=1}^k \int_{\Omega} ( \lambda_i (\mathcal{D}_{\tau}) + c ) d\nu(\tau).
\end{eqnarray}
Hence, $\sum\limits_{i=1}^k \lambda_i (\mathcal{C}) \leq  
\sum\limits_{i=1}^k \int_{\Omega} \lambda_i (\mathcal{D}_{\tau}) d\nu(\tau)$, this is the left statement of Eq.~\eqref{eq1:thm:weak int average thm 4}.
$\hfill \Box$

Next theorem will provide a stronger version of Theorem~\ref{thm:weak int average thm 4} by removing weak majorization conditions. 
\begin{theorem}\label{thm:weak int average thm 5}
Let $\Omega, \nu, \mathcal{C}, \mathcal{D}_\tau$ be defined as the beginning part of Section~\ref{sec:Majorization wtih Integral Average}, and $f: \mathbb{R} \rightarrow [0, \infty)$ be a convex function, we have following two equivalent statements:
\begin{eqnarray}\label{eq1:thm:weak int average thm 5}
\vec{\lambda}(\mathcal{C}) \prec  \int_{\Omega^{m \times p}} \vec{\lambda}(\mathcal{D}_\tau) d\nu^{m \times p}(\tau) \Longleftrightarrow \left\Vert f(\mathcal{C}) \right\Vert_{\rho} \leq 
\int_{\Omega} \left\Vert f(\mathcal{D}_{\tau}) \right\Vert_{\rho}  d\nu(\tau),
\end{eqnarray}
where $\left\Vert \cdot \right\Vert_{\rho}$ is the unitarily invariant norm defined in Eq.~\eqref{eq:def gauge func and general unitarily invariant norm}. 
\end{theorem}
\textbf{Proof:}
We assume that the left statement of Eq.~\eqref{eq1:thm:weak int average thm 5} is true and the function $f$ is a convex function. Again, from Lemma~\ref{lma:Lemma 1 Gen Log Hiai}, we have
\begin{eqnarray}\label{eq2:thm:weak int average thm 5}
\vec{\lambda}(f(\mathcal{A})) = f(\vec{\lambda}(\mathcal{A})) \prec_{w}  f \left( \left(\int_{\Omega^{m \times p}} \vec{\lambda}(\mathcal{D}_\tau) d\nu^{m \times p}(\tau)  \right) \right) \leq \int_{\Omega^{m \times p}} f(\vec{\lambda}(\mathcal{D}_\tau)) d\nu^{m \times p}(\tau),
\end{eqnarray}
then, 
\begin{eqnarray}\label{eq3:thm:weak int average thm 5}
\left\Vert f(\mathcal{A}) \right\Vert_{\rho} &\leq & \rho\left( \int_{\Omega^{m \times p}} f(\vec{\lambda}(\mathcal{D}_\tau)) d\nu^{m \times p}(\tau) \right) \nonumber \\
& \leq & \int_{\Omega}\rho \left( f(\vec{\lambda}(\mathcal{D}_\tau)) \right)d\nu (\tau) = 
\int_{\Omega} \left\Vert f( \mathcal{D}_\tau) \right\Vert_{\rho} d\nu (\tau). 
\end{eqnarray}
This proves the right statement of Eq.~\eqref{eq1:thm:weak int average thm 5}. 

Now, we assume that the right statement of Eq.~\eqref{eq1:thm:weak int average thm 5} is true. From Theorem~\ref{thm:weak int average thm 4}, we already have $\vec{\lambda}(\mathcal{C}) \prec_w  \int_{\Omega^{m \times p}} \vec{\lambda}(\mathcal{D}_\tau) d\nu^{m \times p}(\tau)$.  It is enough to prove $\sum\limits_{i=1}^{m \times p} \lambda_i(\mathcal{C}) \geq \int_{\Omega} \sum\limits_{i=1}^{m \times p} \lambda_i(\mathcal{D}_{\tau}) d \nu(\tau)$. We define a function $f \define \max\{c - x, 0\} $, where $c$ is a positive real constant satisfying $\mathcal{C} \leq c \mathcal{I} $, $\mathcal{D}_{\tau} \leq  c \mathcal{I}$ for all $\tau \in \Omega$ and tensors $c \mathcal{I} - \mathcal{C}, c \mathcal{I} - \mathcal{D}_{\tau}$. If the trace norm is applied, i.e., the sum of the absolute value of all eigenvalues of a symmetric T-product tensor, then the right statement of Eq.~\eqref{eq1:thm:weak int average thm 5} becomes
\begin{eqnarray}\label{eq4:thm:weak int average thm 5}
\sum\limits_{i=1}^{m \times p} \lambda_i \left( c\mathcal{I} - \mathcal{C}\right)  \leq \int_{\Omega} 
\sum\limits_{i=1}^{m \times p} \lambda_i \left( c\mathcal{I} - \mathcal{D}_{\tau}\right) d \nu(\tau).
\end{eqnarray}
The desired inequality  $\sum\limits_{i=1}^{m \times p} \lambda_i(\mathcal{C}) \geq \int_{\Omega} \sum\limits_{i=1}^{m \times p} \lambda_i(\mathcal{D}_{\tau}) d \nu(\tau)$ is established. $\hfill \Box$

\subsection{Log-Majorization with Integral Average}\label{sec:Log-Majorization wtih Integral Average}

The purpose of this section is to consider log-majorization issues for unitarily invariant norms of TPSD T-product tensors. In this section, let $\mathcal{C}, \mathcal{D}_\tau \in \mathbb{R}^{m \times m \times p}$ be TPSD T-product tensors with $m \times p$ nonnegative T-eigenvalues by keeping notations with the same definitions as at the beginning of the Section~\ref{sec:Majorization wtih Integral Average}. For notational convenience, we define the following relation for logarithm vector:
\begin{eqnarray}\label{eq:integral eigen log vector rep}
\left[ \int_{\Omega} \log \lambda_1(\mathcal{D}_\tau) d\nu(\tau), \cdots, \int_{\Omega} \log \lambda_{m \times p}(\mathcal{D}_\tau) d\nu(\tau) \right] \define \int_{\Omega^{m \times p}} \log \vec{\lambda}(\mathcal{D}_\tau) d\nu^{m \times p}(\tau).
\end{eqnarray}

\begin{theorem}\label{thm:weak int average thm 7}
Let $\mathcal{C}, \mathcal{D}_\tau$ be TPSD T-product tensors, $f: (0, \infty) \rightarrow [0,\infty)$ be a continuous function such that the mapping $x \rightarrow \log f(e^x)$ is convex on $\mathbb{R}$, and $g: (0, \infty) \rightarrow [0,\infty)$ be a continuous function such that the mapping $x \rightarrow g(e^x)$ is convex on $\mathbb{R}$ , then  we have following three equivalent statements:
\begin{eqnarray}\label{eq1:thm:weak int average thm 7}
\vec{\lambda}(\mathcal{C}) &\prec_{w \log}& \exp  \int_{\Omega^{m \times p}} \log \vec{\lambda}(\mathcal{D}_\tau) d\nu^{m \times p}(\tau);
\end{eqnarray}
\begin{eqnarray}\label{eq2:thm:weak int average thm 7}
\left\Vert f(\mathcal{C}) \right\Vert_{\rho} &\leq &
\exp \int_{\Omega} \log \left\Vert f(\mathcal{D}_{\tau}) \right\Vert_{\rho}  d\nu(\tau);
\end{eqnarray}
\begin{eqnarray}\label{eq3:thm:weak int average thm 7}
\left\Vert g(\mathcal{C}) \right\Vert_{\rho} &\leq &
\int_{\Omega} \left\Vert g(\mathcal{D}_{\tau}) \right\Vert_{\rho}  d\nu(\tau).
\end{eqnarray}
\end{theorem}
\textbf{Proof:}
The roadmap of this proof is to prove equivalent statements between Eq.~\eqref{eq1:thm:weak int average thm 7} and Eq.~\eqref{eq2:thm:weak int average thm 7} first, followed by equivalent statements between Eq.~\eqref{eq1:thm:weak int average thm 7} and Eq.~\eqref{eq3:thm:weak int average thm 7}. 

\textbf{Eq.~\eqref{eq1:thm:weak int average thm 7} $\Longrightarrow$ Eq.~\eqref{eq2:thm:weak int average thm 7}}

There are two cases to be discussed in this part of proof: $\mathcal{C}, \mathcal{D}_\tau$ are TPD tensors, and $\mathcal{C}, \mathcal{D}_\tau$ are TPSD T-product tensors. At the beginning, we consider the case that $\mathcal{C}, \mathcal{D}_\tau$ are TPD tensors.

Since $\mathcal{D}_\tau$ are positive, we can find $\varepsilon > 0$ such that $\mathcal{D}_{\tau} \geq \varepsilon \mathcal{I}$ for all $\tau \in \Omega$. From Eq.~\eqref{eq1:thm:weak int average thm 7}, the convexity of $\log f (e^x)$ and Lemma~\ref{lma:Lemma 1 Gen Log Hiai}, we have 
\begin{eqnarray}
\vec{\lambda} \left( f ( \mathcal{C})\right)  = f \left(\exp \left( \log \vec{\lambda} (\mathcal{C}) \right) \right) &\prec _w &  f \left(\exp    \int_{\Omega^{m \times p}} \vec{\lambda}(\mathcal{D}_\tau) d\nu^{m \times p}(\tau)             \right) \nonumber \\
& \leq &  \exp \left( \int_{\Omega^{m \times p}} \log f \left( \vec{\lambda}(\mathcal{D}_\tau) \right)  d\nu^{m \times p}(\tau) \right).
\end{eqnarray}
Then, from Eq.~\eqref{eq:def gauge func and general unitarily invariant norm}, we obtain
\begin{eqnarray}\label{eq4:thm:weak int average thm 7}
\left\Vert f (\mathcal{C}) \right\Vert _{\rho}
& \leq &  \rho\left(\exp \left( \int_{\Omega^{m \times p}} \log  f \left( \vec{\lambda}(\mathcal{D}_\tau) \right) d\nu^{m \times p}(\tau) \right)   \right).
\end{eqnarray}


From the function $f$ properties, we can assume that $f(x) > 0$ for any $x > 0$. Then, we have 
following bounded and continous maps on $\Omega$: $\tau \rightarrow \log f (\lambda_i (\mathcal{D}_{\tau}))  $ for $i \in \{1,2,\cdots, m \times p \}$, and $\tau \rightarrow \log \left\Vert f (\mathcal{D}_{\tau}) \right\Vert_{\rho}$. Because we have $\nu (\Omega) = 1$ and $\sigma$-compactness of $\Omega$, we have $\tau_{k}^{(n)} \in \Omega$ and $\alpha_{k}^{(n)}$ for $k \in \{1,2,\cdots, n\}$ and $n \in \mathbb{N}$ with $\sum\limits_{k=1}^{n} \alpha_{k}^{(n)} = 1$ such that 
\begin{eqnarray}\label{eq:35}
\int_{\Omega} \log f (\lambda_i ( \mathcal{D}_{\tau} ) ) d \nu (\tau) = \lim\limits_{n \rightarrow \infty} \sum\limits_{k=1}^{n} \alpha_{k}^{(n)}  \log f (\lambda_i (\mathcal{D}_{\tau_k^{(n)} }))   , \mbox{for $i \in \{1,2,\cdots  m \times p \}$};
\end{eqnarray}
and 
\begin{eqnarray}\label{eq:36}
\int_{\Omega} \log \left\Vert f (\mathcal{D}_{\tau}) \right\Vert_{\rho} d \nu (\tau) = \lim\limits_{n \rightarrow \infty} \sum\limits_{k=1}^{n} \alpha_{k}^{(n)}  \log \left\Vert f (\mathcal{D}_{\tau_k^{(n)}    }) \right\Vert_{\rho} .
\end{eqnarray}
By taking the exponential at both sides of Eq.~\eqref{eq:35} and apply the gauge function $\rho$, we have
\begin{eqnarray}\label{eq:37}
\rho \left( \exp \int_{\Omega^{m \times p}} \log f (\vec{\lambda} ( \mathcal{D}_{\tau} ) ) d \nu^{m \times p} (\tau)  \right)= \lim\limits_{n \rightarrow \infty} \rho\left(  \prod\limits_{k=1}^{n}  f  \left( \vec{\lambda} \left(\mathcal{D}_{\tau_k^{(n)} } \right)  \right)^{ \alpha_{k}^{(n)} }  \right).
\end{eqnarray}
Similarly, by taking the exponential at both sides of Eq.~\eqref{eq:36}, we have
\begin{eqnarray}\label{eq:38}
\exp \left( \int_{\Omega} \log \left\Vert f (\mathcal{D}_{\tau}) \right\Vert_{\rho} d \nu (\tau) \right) = \lim\limits_{n \rightarrow \infty} \prod \limits_{k=1}^{n} \left\Vert f \left( \mathcal{D}_{\tau_k^{(n)}    } \right) \right\Vert^{\alpha_{k}^{(n)}}_{\rho} .
\end{eqnarray}
From Lemma~\ref{lma:Holder inquality for gauge function}, we have 
\begin{eqnarray}\label{eq:41}
\rho \left(  \prod\limits_{k=1}^{n}  f  \left( \vec{\lambda} \left(\mathcal{D}_{\tau_k^{(n)} } \right)  \right)^{ \alpha_{k}^{(n)} }  \right) & \leq & \prod \limits_{k=1}^{n} \rho^{  \alpha_{k}^{(n)}    } \left( f \left( \vec{\lambda} \left( \mathcal{D}_{ \tau_{k}^{(n)}  }\right) 
 \right) \right) \nonumber \\
&=& \prod \limits_{k=1}^{n} \rho^{  \alpha_{k}^{(n)}    } \left( \vec{\lambda}  \left( f  \left( \mathcal{D}_{ \tau_{k}^{(n)}  }\right) 
 \right) \right) \nonumber \\
&=& \prod \limits_{k=1}^{n} \left\Vert f  \left( \mathcal{D}_{ \tau_{k}^{(n)}  } \right) \right\Vert_{\rho}^{\alpha_{k}^{(n)}}
\end{eqnarray}

From Eqs.~\eqref{eq:37},~\eqref{eq:38} and~\eqref{eq:41}, we have 
\begin{eqnarray}\label{eq:42}
\rho \left( \exp \int_{\Omega^{m \times p}} \log f (\vec{\lambda} ( \mathcal{D}_{\tau} ) ) d \nu^{m \times p} (\tau)  \right) \leq \exp \int_{\Omega} \log \left\Vert f(\mathcal{D}_{\tau})\right\Vert_{\rho} d \nu(\tau).
\end{eqnarray}
Then, Eq.~\eqref{eq2:thm:weak int average thm 7} is proved from Eqs.~\eqref{eq4:thm:weak int average thm 7} and~\eqref{eq:42}.

Next, we consider that $\mathcal{C}, \mathcal{D}_\tau$ are TPSD T-product tensors. For any $\delta > 0$, we have following log-majorization relation:
\begin{eqnarray}
\prod\limits_{i=1}^k \left( \lambda_i (\mathcal{C}) + \epsilon_{\delta} \right) 
&\leq& \prod\limits_{i=1}^k \exp  \int_{\Omega} \log \left( \lambda_i(\mathcal{D}_\tau) + \delta\right) d \nu (\tau),
\end{eqnarray}
where $\epsilon_{\delta} > 0$ and $k \in \{1,2,\cdots r \}$. Then, we can apply the previous case result about TPD tensors to TPD tensors $\mathcal{C} + \epsilon_{\delta} \mathcal{I}$ and $\mathcal{D}_\tau + \delta \mathcal{I}$, and get 
\begin{eqnarray}\label{eq:46}
\left\Vert f (\mathcal{C}) + \epsilon_{\delta}  \mathcal{I} \right\Vert_{\rho} 
&\leq& \exp \int_{\Omega} \log \left\Vert f (\mathcal{D}_{\tau}) + \delta  \mathcal{I} \right\Vert_{\rho} 
d \nu (\tau)
\end{eqnarray}
As $\delta \rightarrow 0$, Eq.~\eqref{eq:46} will give us Eq.~\eqref{eq2:thm:weak int average thm 7} for TPSD T-product tensors.  

\textbf{Eq.~\eqref{eq1:thm:weak int average thm 7} $\Longleftarrow$ Eq.~\eqref{eq2:thm:weak int average thm 7}}

We consider TPD tensors at first phase by assuming that $\mathcal{D}_{\tau}$ are 
TPD T-product tensors for all $\tau \in \Omega$. We may also assume that the tensor $\mathcal{C}$ is a TPD T-product tensor. Since if this is a TPSD T-product tensor, i.e., some $\lambda_i = 0$, we always have following inequality valid:
\begin{eqnarray}
\prod\limits_{i=1}^k \lambda_i (\mathcal{C}) \leq \prod\limits_{i=1}^k 
\exp \int_{\Omega} \log \lambda_i (\mathcal{D}_{\tau}) d \nu (\tau)
\end{eqnarray}

If we apply $f(x) = x^p$ for $p > 0$ and $\left\Vert \cdot \right\Vert_{\rho}$ as Ky Fan $k$-norm in Eq.~\eqref{eq2:thm:weak int average thm 7}, we have 
\begin{eqnarray}\label{eq:50}
\log \sum\limits_{i=1}^k \lambda^p_i \left(\mathcal{C}\right) \leq \int_{\Omega} \log \sum\limits_{i=1}^k \lambda_i^p\left( \mathcal{D}_{\tau} \right) d \nu (\tau).
\end{eqnarray}
If we add $\log \frac{1}{k}$ and multiply $\frac{1}{p}$ at both sides of Eq.~\eqref{eq:50}, we have 
\begin{eqnarray}\label{eq:51}
\frac{1}{p}\log \left( \frac{1}{k}\sum\limits_{i=1}^k \lambda^p_i \left(\mathcal{C}\right) \right)\leq \int_{\Omega} \frac{1}{p} \log \left( \frac{1}{k}\sum\limits_{i=1}^k \lambda_i^p\left( \mathcal{D}_{\tau} \right) \right) d \nu (\tau).
\end{eqnarray}
From L'Hopital's Rule, if $p \rightarrow 0$, we have 
\begin{eqnarray}\label{eq:52}
\frac{1}{p}\log \left( \frac{1}{k}\sum\limits_{i=1}^k \lambda^p_i \left(\mathcal{C}\right) \right) \rightarrow \frac{1}{k} \sum\limits_{i=1}^k \log \lambda_i (\mathcal{C}),
\end{eqnarray}
and 
\begin{eqnarray}\label{eq:53}
\frac{1}{p}\log \left( \frac{1}{k}\sum\limits_{i=1}^k \lambda^p_i \left(\mathcal{D}_{\tau}\right) \right) \rightarrow \frac{1}{k} \sum\limits_{i=1}^k \log \lambda_i (\mathcal{D}_{\tau}),
\end{eqnarray}
where $\tau \in \Omega$. Appling Eqs.~\eqref{eq:52} and~\eqref{eq:53} into Eq.~\eqref{eq:51} and taking $p \rightarrow 0$, we have 
\begin{eqnarray}
\sum\limits_{i=1}^k \lambda_i (\mathcal{C}) \leq \int_{\Omega} \sum\limits_{i=1}^k 
 \log \lambda_i (\mathcal{D}_{\tau}) d \nu (\tau).
\end{eqnarray}
Therefore, Eq.~\eqref{eq1:thm:weak int average thm 7} is true for TPD tensors. 

For TPSD T-product tensors $\mathcal{D}_{\tau}$, since Eq.~\eqref{eq2:thm:weak int average thm 7} is valid for $\mathcal{D}_{\tau} + \delta \mathcal{I}$ for any $\delta > 0$, we can apply the previous case result about TPD tensors to $\mathcal{D}_{\tau} + \delta \mathcal{I}$ and obtain
\begin{eqnarray}
\prod\limits_{i=1}^k \lambda_i (\mathcal{C}) \leq \prod\limits_{i=1}^k  \exp 
\int_{\Omega}   \log  \left( \lambda_i (\mathcal{D}_{\tau}) + \delta \right) d \nu (\tau),
\end{eqnarray}
where $k \in \{1,2,\cdots, r\}$. Eq.~\eqref{eq1:thm:weak int average thm 7} is still true for TPSD T-product tensors as $\delta \rightarrow 0$.

\textbf{Eq.~\eqref{eq1:thm:weak int average thm 7} $\Longrightarrow$ Eq.~\eqref{eq3:thm:weak int average thm 7}}

If $\mathcal{C}, \mathcal{D}_\tau$ are TPD tensors, and $\mathcal{D}_{\tau} \geq \delta \mathcal{I}$ for all $\tau \in \Omega$. From Eq.~\eqref{eq1:thm:weak int average thm 7}, we have 
\begin{eqnarray}
\vec{\lambda} (\log \mathcal{C}) = \log \vec{\lambda}(\mathcal{C}) \prec_{w}
\int_{\Omega^{m \times p}} \log \vec{\lambda}(\mathcal{D}_{\tau}) d \nu^{m \times p} (\tau) = 
\int_{\Omega^{m \times p}} \vec{\lambda}( \log \mathcal{D}_{\tau}) d \nu^{m \times p} (\tau).
\end{eqnarray}
If we apply Theorem~\ref{thm:weak int average thm 4} to $\log \mathcal{C}$, $\log \mathcal{D}_{\tau}$ with function $f(x) = g(e^x)$, where $g$ is used in Eq.~\eqref{eq3:thm:weak int average thm 7}, Eq.~\eqref{eq3:thm:weak int average thm 7} is implied. 

If $\mathcal{C}, \mathcal{D}_\tau$ are TPSD T-product tensors and any $\delta > 0$, we can find $\epsilon_{\delta} \in (0, \delta)$ to satisfy following:
\begin{eqnarray}\label{eq:45}
\prod\limits_{i=1}^k\left(\lambda_i(\mathcal{C}) + \epsilon_{\delta}\right) \leq 
\prod\limits_{i=1}^k \exp \int_{\Omega}   \log \left( \lambda_i(\mathcal{D}_{\tau}) + \delta  \right) d \nu (\tau).
\end{eqnarray}
Then, from TPD T-product tensor case, we have 
\begin{eqnarray}\label{eq:45-1}
\left\Vert g( \mathcal{C} + \epsilon_{\delta} \mathcal{I} ) \right\Vert_{\rho}
\leq \int_{\Omega} \left\Vert   g( \mathcal{D}_{\tau} + \delta \mathcal{I} )    \right\Vert_{\rho}
d \nu (\tau).
\end{eqnarray}
Eq.~\eqref{eq3:thm:weak int average thm 7} is obtained by taking $\delta \rightarrow 0$ in Eq.~\eqref{eq:45-1}. 

\textbf{Eq.~\eqref{eq1:thm:weak int average thm 7} $\Longleftarrow$ Eq.~\eqref{eq3:thm:weak int average thm 7}}

For $k \in \{1,2,\cdots, r \}$, if we apply $g(x) = \log (\delta + x )$, where $\delta >0$, and Ky Fan $k$-norm in Eq.~\eqref{eq3:thm:weak int average thm 7}, we have 
\begin{eqnarray}
\sum\limits_{i=1}^k \log \left(\delta + \lambda_i \left(\mathcal{C} \right) \right)
\leq \sum\limits_{i=1}^k \int_{\Omega} \log \left( \delta + \lambda_{i}(\mathcal{D}_{\tau}) \right) d \nu (\tau).
\end{eqnarray}
Then, we have following relation as $\delta \rightarrow 0$:
\begin{eqnarray}
\sum\limits_{i=1}^k \log \lambda_i \left(\mathcal{C} \right) 
\leq \sum\limits_{i=1}^k \int_{\Omega} \log  \lambda_{i}(\mathcal{D}_{\tau}) d \nu (\tau).
\end{eqnarray}
Therefore, Eq.~\eqref{eq1:thm:weak int average thm 7} ccan be derived from Eq.~\eqref{eq3:thm:weak int average thm 7}. $\hfill \Box$

Next theorem will extend Theorem~\ref{thm:weak int average thm 7} to non-weak version.

\begin{theorem}\label{thm:int log average thm 10}
Let $\mathcal{C}, \mathcal{D}_\tau$ be TPSD T-product tensors with $\int_{\Omega} \left\Vert \mathcal{D}_{\tau}^{-p}\right\Vert_\rho d \nu (\tau) < \infty$ for any $p > 0$, $f: (0, \infty) \rightarrow [0,\infty)$ be a continuous function such that the mapping $x \rightarrow \log f(e^x)$ is convex on $\mathbb{R}$, and $g: (0, \infty) \rightarrow [0,\infty)$ be a continuous function such that the mapping $x \rightarrow g(e^x)$ is convex on $\mathbb{R}$ , then  we have following three equivalent statements:
\begin{eqnarray}\label{eq1:thm:int average thm 10}
\vec{\lambda}(\mathcal{C}) &\prec_{\log}& \exp  \int_{\Omega^{m \times p}} \log \vec{\lambda}(\mathcal{D}_\tau) d\nu^{m \times p}(\tau);
\end{eqnarray}
\begin{eqnarray}\label{eq2:thm:int average thm 10}
\left\Vert f(\mathcal{C}) \right\Vert_{\rho} &\leq &
\exp \int_{\Omega} \log \left\Vert f(\mathcal{D}_{\tau}) \right\Vert_{\rho}  d\nu(\tau);
\end{eqnarray}
\begin{eqnarray}\label{eq3:thm:int average thm 10}
\left\Vert g(\mathcal{C}) \right\Vert_{\rho} &\leq &
\int_{\Omega} \left\Vert g(\mathcal{D}_{\tau}) \right\Vert_{\rho}  d\nu(\tau).
\end{eqnarray}
\end{theorem}
\textbf{Proof:}

The proof plan is similar to the proof in Theorem~\ref{thm:weak int average thm 7}. We prove the equivalence between Eq.~\eqref{eq1:thm:int average thm 10} and Eq.~\eqref{eq2:thm:int average thm 10} first, then prove the equivalence between Eq.~\eqref{eq1:thm:int average thm 10} and Eq.~\eqref{eq3:thm:int average thm 10}. 

\textbf{Eq.~\eqref{eq1:thm:int average thm 10} $\Longrightarrow$ Eq.~\eqref{eq2:thm:int average thm 10}}

First, we assume that $\mathcal{C}, \mathcal{D}_\tau$ are TPD T-product tensors with $\mathcal{D}_\tau \geq \delta \mathcal{I}$ for all $\tau \in \Omega$. The corresponding part of the proof in Theorem~\ref{thm:weak int average thm 7} about TPD tensors $\mathcal{C}, \mathcal{D}_\tau$ can be applied here. 

For case that $\mathcal{C}, \mathcal{D}_\tau$ are TPSD T-product tensors, we have 
\begin{eqnarray}
\prod\limits_{i=1}^k \lambda_i (\mathcal{C}) \leq \prod\limits_{i=1}^k 
\exp \int_{\Omega} \log \left( \lambda_i (\mathcal{D}_{\tau}) + \delta_n \right)d \nu (\tau), 
\end{eqnarray}
where $\delta_n > 0$ and $\delta_n \rightarrow 0$. Because $\int_{\Omega^{m \times p}} \log \left( \vec{\lambda} (\mathcal{D}_{\tau}) + \delta_n \right) d \nu^{m \times p} (\tau) \rightarrow \int_{\Omega^{m \times p}} \log \vec{\lambda} (\mathcal{D}_\tau)  d \nu^{m \times p} (\tau) $ as $n \rightarrow \infty$, from Lemma~\ref{lma:Lemma 12 Gen Log Hiai}, we can find $\mathbf{a}^{(n)}$ with $n \geq n_0$ such that $a^{(n)}_1 \geq \cdots \geq a^{(n)}_r > 0$, $\mathbf{a}^{(n)} \rightarrow \vec{\lambda}(\mathcal{C})$ and $ \mathbf{a}^{(n)} \prec_{\log}  \exp \int_{\Omega^{m \times p}} \log \vec{\lambda} \left(\mathcal{D}_{\tau} + \delta_n \mathcal{I} \right) d \nu^{m \times p} (\tau)$

Selecting $\mathcal{C}^{(n)}$ with $\vec{\lambda} ( \mathcal{C}^{(n)})  = \mathbf{a}^{(n)} $ and applying TPD tensors case to $\mathcal{C}^{(n)}$ and $\mathcal{D}_{\tau} + \delta_n \mathcal{I}$, we obtain
\begin{eqnarray}\label{eq:74}
\left\Vert f (\mathcal{C}^{(n)}) \right\Vert_{\rho} \leq \exp \int_{\Omega} \log \left\Vert f (\mathcal{D}_{\tau} + \delta_n \mathcal{I}) \right\Vert_{\rho} d \nu (\tau)
\end{eqnarray}
where $n \geq n_0$.

There are two situations for the function $f$ near $0$: $f(0^{+}) < \infty$ and $f(0^{+}) = \infty$. For the case with $f(0^{+}) < \infty$, we have 
\begin{eqnarray}\label{eq:75-1}
\left\Vert f (\mathcal{C}^{(n)} )\right\Vert_{\rho} = \rho( f (\mathbf{a}^{(n)}))
\rightarrow \rho (f ( \vec{\lambda}(\mathcal{C}))) = \left\Vert f (\mathcal{C})\right\Vert_{\rho}, 
\end{eqnarray}
and
\begin{eqnarray}\label{eq:75-2}
\left\Vert f (\mathcal{D}_{\tau} + \delta_n \mathcal{I} )\right\Vert_{\rho} 
\rightarrow \left\Vert f (\mathcal{D}_{\tau})\right\Vert_{\rho}, 
\end{eqnarray}
where $\tau \in \Omega$ and $n \rightarrow \infty$. From Fatou–Lebesgue theorem, we then have 
\begin{eqnarray}\label{eq:76}
\limsup\limits_{n \rightarrow \infty} \int_{\Omega} \log \left\Vert f (\mathcal{D}_{\tau} + \delta_n \mathcal{I} )\right\Vert_{\rho} d \nu (\tau) \leq \int_{\Omega} \log \left\Vert f(\mathcal{D}_{\tau}) \right\Vert_{\rho}.
\end{eqnarray}
By taking $n \rightarrow \infty$ in Eq.~\eqref{eq:74} and using Eqs.~\eqref{eq:75-1},~\eqref{eq:75-2},~\eqref{eq:76}, we have Eq.~\eqref{eq2:thm:int average thm 10} for case that $f(0^{+}) < \infty$.

For the case with $f(0^{+}) = \infty$, we assume that $\int_{\Omega} \log \left\Vert f (\mathcal{D}_{\tau}) \right\Vert_{\rho} d \nu (\tau) < \infty$ (since the inequality in Eq.~\eqref{eq2:thm:int average thm 10} is always true for $\int_{\Omega} \log \left\Vert f (\mathcal{D}_{\tau}) \right\Vert_{\rho} d \nu (\tau) = \infty$). Since $f$ is decreasing on $(0, \epsilon)$ for some $\epsilon > 0$. We claim that the following relation is valid: there are two constants $a, b > 0$ such that 
\begin{eqnarray}\label{eq:77}
a \leq \left\Vert f (\mathcal{D}_{\tau} + \delta_n \mathcal{I}) \right\Vert_{\rho}
\leq \left\Vert f (\mathcal{D}_{\tau}) \right\Vert_{\rho} + b,
\end{eqnarray}
for all $\tau \in \Omega$ and $n \geq n_0$. If Eq.~\eqref{eq:77} is valid and $\int_{\Omega} \log \left\Vert f (\mathcal{D}_{\tau}) \right\Vert_{\rho} d \nu (\tau) < \infty$, from Lebesgue's dominated convergence theorem, we also have Eq.~\eqref{eq2:thm:int average thm 10} for case that $f(0^{+}) = \infty$ by taking $n \rightarrow \infty$ in Eq.~\eqref{eq:74}. 

Below, we will prove the claim stated by Eq.~\eqref{eq:77}. By the uniform boundedness of tensors $\mathcal{D}_{\tau}$, there is a constant $\kappa >0$ such that 
\begin{eqnarray}
0 < \mathcal{D}_{\tau} + \delta_n \mathcal{I} \leq \kappa \mathcal{I},
\end{eqnarray}
where $\tau \in \Omega$ and $ n \geq n_0$. We may assume that $\mathcal{D}_\tau$ is TPD tensors because $\left\Vert f (\mathcal{D}_{\tau}) \right\Vert_{\rho} = \infty$, i.e., Eq.~\eqref{eq:77} being true automatically, when $\mathcal{D}_\tau$ is TPSD T-product tensors. From Theorem~\ref{thm:T eigenvalue decomp}, we have 
\begin{eqnarray}
f(\mathcal{D}_{\tau} + \delta_n \mathcal{I}) &=&  \sum\limits_{i', \mbox{s.t. $\lambda_{i'}(\mathcal{D}_{\tau})+  \delta_n < \epsilon$}} f(\lambda_{i'}(\mathcal{D}_{\tau}) + \delta_n ) \mathcal{U}_{i'} \star \mathcal{U}^{H}_{i'} + \nonumber \\
&  & \sum\limits_{j', \mbox{s.t. $\lambda_{j'}(\mathcal{D}_{\tau}) +  \delta_n 
\geq \epsilon$}} f(\lambda_{j'}(\mathcal{D}_{\tau}) + \delta_n ) \mathcal{U}_{j'} \star \mathcal{U}^{H}_{j'}  \nonumber \\
&\leq &  \sum\limits_{i', \mbox{s.t. $\lambda_{i'}(\mathcal{D}_{\tau}) +  \delta_n < \epsilon$}} f(\lambda_{i'}(\mathcal{D}_{\tau}) ) \mathcal{U}_{i'} \star\mathcal{U}^{H}_{i'} + \nonumber \\
&  & \sum\limits_{j', \mbox{s.t. $\lambda_{j'}(\mathcal{D}_{\tau})  +  \delta_n 
\geq \epsilon$}} f(\lambda_{j'}(\mathcal{D}_{\tau}) + \delta_n ) \mathcal{U}_{j'} \star \mathcal{U}^{H}_{j'}  \nonumber \\
&\leq & f(\mathcal{D}_{\tau}) + \sum\limits_{j', \mbox{s.t. $\lambda_{j'}(\mathcal{D}_{\tau}) +  \delta_n \geq \epsilon$}} f(\lambda_{j'}(\mathcal{D}_{\tau}) + \delta_n ) \mathcal{U}_{j'} \star \mathcal{U}^{H}_{j'}.
\end{eqnarray}
Therefore, the claim in Eq.~\eqref{eq:77} follows by the triangle inequality for $\left\Vert \cdot \right\Vert_{\rho}$ and $f(\lambda_{j'}(\mathcal{D}_{\tau}) + \delta_n )  < \infty$ for $\lambda_{j'}(\mathcal{D}_{\tau}) +  \delta_n \geq \epsilon$. 

\textbf{Eq.~\eqref{eq1:thm:int average thm 10} $\Longleftarrow$ Eq.~\eqref{eq2:thm:int average thm 10}}

The weak majorization relation 
\begin{eqnarray}\label{eq:82}
\prod\limits_{i=1}^{k} \lambda_i (\mathcal{C}) \leq \prod\limits_{i=1}^{k} \exp \int_{\Omega} \log \lambda_i (\mathcal{D}_\tau) d \nu (\tau),
\end{eqnarray}
is valid for $k < m \times p$ from Eq.~\eqref{eq1:thm:weak int average thm 7} $\Longrightarrow$ Eq.~\eqref{eq2:thm:weak int average thm 7} in Theorem~\ref{thm:weak int average thm 7}.  We wish to prove that Eq.~\eqref{eq:82} becomes equal for $k =  m \times p$. It is equivalent to prove that 
\begin{eqnarray}\label{eq:83}
\log \det (\mathcal{C}) \geq   \int_{\Omega} \log \det (\mathcal{D}_{\tau}) d \nu (\tau),
\end{eqnarray}
where $\det ( \cdot )$ is defined by Eq.~\eqref{eq:def T prod tensor determinant}. We can assume that $ \int_{\Omega} \log \det (\mathcal{D}_{\tau}) d \nu (\tau) \geq - \infty$ since Eq.~\eqref{eq:83} is true for  $ \int_{\Omega} \log \det (\mathcal{D}_{\tau}) d \nu (\tau) = - \infty$. Then, $\mathcal{D}_\tau$ are TPD tensors. 

If we scale tensors $\mathcal{C}, \mathcal{D}_{\tau}$ as $a \mathcal{C}, a\mathcal{D}_{\tau}$ by some $a >0$, we can assume $\mathcal{D}_{\tau} \leq \mathcal{I}$ and $\lambda_i(\mathcal{D}_\tau) \leq 1$ for all $ \tau \in \Omega$ and $ i \in \{1,2,\cdots, m \times p\}$. Then for any $p >0$, we have 
\begin{eqnarray}
\frac{1}{m \times p} \left\Vert \mathcal{D}_{\tau}^{-\varrho} \right\Vert_1 \leq \lambda^{-\varrho}_r (\mathcal{D}_{\tau} ) \leq ( \det  (\mathcal{D}_{\tau})  )^{-\varrho},
\end{eqnarray}
and 
\begin{eqnarray}\label{eq:85}
\frac{1}{\varrho} \log \left( \frac{\left\Vert \mathcal{D}^{-\varrho}_\tau \right\Vert_1 }{m \times p}\right)
\leq - \log \det  (\mathcal{D}_{\tau}). 
\end{eqnarray}
If we use tensor trace norm, represented by $\left\Vert \cdot \right\Vert_1$,  as unitarily invariant tensor norm and $f(x) = x^{-\varrho}$ for any $\varrho > 0$ in Eq.~\eqref{eq2:thm:int average thm 10}, we obtain
\begin{eqnarray}\label{eq:86}
\log \left\Vert \mathcal{C}^{-\varrho} \right\Vert_1 \leq \int_{\Omega} \log \left\Vert \mathcal{D}^{-\varrho}_{\tau} \right\Vert_1 d \nu(\tau).
\end{eqnarray}
By adding $\log \frac{1}{m \times p}$ and multiplying $\frac{1}{\varrho}$ for both sides of Eq.~\eqref{eq:86}, we have 
\begin{eqnarray}\label{eq:87}
\frac{1}{\varrho}\log \left( \frac{\left\Vert \mathcal{C}^{-\varrho} \right\Vert_1 }{m \times p} \right)
\leq \int_{\Omega}\frac{1}{\varrho} \log \left(  \frac{\left\Vert \mathcal{D}_\tau^{-\varrho} \right\Vert_1 }{m \times p}           \right) d \nu (\tau)
\end{eqnarray}
Similar to Eqs.~\eqref{eq:52} and~\eqref{eq:53}, we have following two relations as $\varrho \rightarrow 0$:
\begin{eqnarray}\label{eq:88}
\frac{1}{\varrho}\log \left( \frac{\left\Vert \mathcal{C}^{-\varrho} \right\Vert_1 }{m \times p} \right) \rightarrow \frac{- 1}{m \times p} \log \det (\mathcal{C}),
\end{eqnarray}
and 
\begin{eqnarray}\label{eq:89}
\frac{1}{\varrho}\log \left( \frac{\left\Vert \mathcal{D}_{\tau}^{-\varrho} \right\Vert_1}{m \times p}  \right) \rightarrow \frac{- 1}{m \times p} \log \det (\mathcal{D}_\tau).
\end{eqnarray}
From Eq.~\eqref{eq:85} and Lebesgue's dominated convergence theorem, we have 
\begin{eqnarray}\label{eq:90}
\lim\limits_{\varrho  \rightarrow 0} \int_{\Omega}\frac{1}{\varrho}\log \left( \frac{\left\Vert \mathcal{D}_{\tau}^{-\varrho} \right\Vert_1}{m \times p}  \right) d \nu (\tau)= \frac{-1}{m \times p} \int_{\Omega} \log \det(\mathcal{D}_{\tau})     \nu (\tau) 
\end{eqnarray}
Finally, we have Eq.~\eqref{eq:83} from Eqs.~\eqref{eq:87} and~\eqref{eq:90}.

\textbf{Eq.~\eqref{eq1:thm:int average thm 10} $\Longrightarrow$ Eq.~\eqref{eq3:thm:int average thm 10}}

First, we assume that $\mathcal{C}, \mathcal{D}_\tau$ are TPD tensors and $\mathcal{D}_{\tau} \geq \delta \mathcal{I}$ for $\tau \in \Omega$. From Eq.~\eqref{eq1:thm:int average thm 10}, we can apply Theorem~\ref{thm:weak int average thm 5} to $\log \mathcal{C}, \log \mathcal{D}_\tau$ and $f(x) = g(e^x)$ to obtain Eq.~\eqref{eq3:thm:int average thm 10}. 

For $\mathcal{C}, \mathcal{D}_\tau$ are TPSD T-product tensors, we can choose  $\mathbf{a}^{(n)}$ and corresponding  $\mathcal{C}^{(n)}$ for $n \geq n_0$ given $\delta_n \rightarrow 0$ with $\delta_n > 0$ as the proof in Eq.~\eqref{eq1:thm:int average thm 10} $\Longrightarrow$ Eq.~\eqref{eq2:thm:int average thm 10}. Since tensors $\mathcal{C}^{(n)}, \mathcal{D}_\tau + \delta_n \mathcal{I}$ are TPD T-product tensors, we then have 
\begin{eqnarray}\label{eq:92}
\left \Vert g ( \mathcal{C}^{(n)}  ) \right\Vert_{\rho} \leq \int_{\Omega} \left\Vert g ( \mathcal{D}_\tau + \delta_n \mathcal{I} ) \right\Vert_{\rho} d \nu (\tau).
\end{eqnarray}
If $g(0^+) < \infty$, Eq.~\eqref{eq3:thm:int average thm 10} is obtained from Eq.~\eqref{eq:92} by taking $n \rightarrow \infty$. On the other hand, if $g(0^+) = \infty$, we can apply the argument similar to the portion about $f(0^+) = \infty$ in the proof for Eq.~\eqref{eq1:thm:int average thm 10} $\Longrightarrow$ Eq.~\eqref{eq2:thm:int average thm 10} to get $a, b > 0$ such that 
\begin{eqnarray}\label{eq:92 infty}
a \leq \left\Vert g ( \mathcal{D}_\tau + \delta_n \mathcal{I} ) \right\Vert_\rho \leq   \left\Vert  g ( \mathcal{D}_\tau ) \right\Vert_\rho + b,
\end{eqnarray}
for all $\tau \in \Omega$ and $n \geq n_0$. Since the case that $\int_{\Omega} \left\Vert  g ( \mathcal{D}_\tau ) \right\Vert_{\rho} d \nu (\tau) = \infty$ will have Eq.~\eqref{eq3:thm:int average thm 10}, we only consider the case that $\int_{\Omega} \left\Vert  g ( \mathcal{D}_\tau ) \right\Vert_{\rho} d \nu (\tau) < \infty$. Then, we have Eq.~\eqref{eq3:thm:int average thm 10} from Eqs.~\eqref{eq:92},~\eqref{eq:92 infty} and Lebesgue's dominated convergence theorem.

\textbf{Eq.~\eqref{eq1:thm:int average thm 10} $\Longleftarrow$ Eq.~\eqref{eq3:thm:int average thm 10}}

The weak majorization relation
\begin{eqnarray}\label{eq:94}
\sum\limits_{i=1}^{k} \log \lambda_i (\mathcal{C}) \leq 
\sum\limits_{i=1}^{k} \int_{\Omega} \log \lambda_i (\mathcal{D}_\tau) d \nu (\tau)
\end{eqnarray}
is true from the implication from Eq.~\eqref{eq1:thm:weak int average thm 7} to Eq.~\eqref{eq3:thm:weak int average thm 7} in Theorem~\ref{thm:weak int average thm 7}. We have to show that this relation becomes identity for $k=m \times p$. If we apply $\left\Vert \cdot  \right\Vert_{\rho} = \left\Vert \cdot \right\Vert_1$ and $g(x) = x^{-\varrho}$ for any $\varrho > 0$ in Eq.~\eqref{eq3:thm:int average thm 10}, we have 
\begin{eqnarray}\label{eq:95}
\frac{1}{\varrho} \log \left( \frac{ \left\Vert \mathcal{C}^{-\varrho}\right\Vert_1  }{m \times p} \right)
\leq  \frac{1}{\varrho} \log \left( \int_{\Omega} \frac{\left\Vert \mathcal{D}_\tau^{-\varrho}\right\Vert_1}{m \times p}   d \nu(\tau)  \right).
\end{eqnarray}
Then, we will get 
\begin{eqnarray}\label{eq:96}
\frac{- \log \det (\mathcal{C})}{m \times p} &=& \lim\limits_{\varrho \rightarrow 0} \frac{1}{\varrho} \log \left( \frac{ \left\Vert \mathcal{C}^{-\varrho}\right\Vert_1  }{m \times p} \right) \nonumber \\
& \leq & \lim\limits_{\varrho \rightarrow 0} \frac{1}{p} \log \left( \int_{\Omega} \frac{\left\Vert \mathcal{D}_\tau^{-\varrho}\right\Vert_1}{m \times p}   d \nu(\tau)  \right)=_1 \frac{ - \int_{\Omega} \log \det (\mathcal{D}_{\tau}) d \nu (\tau)  }{m \times p},
\end{eqnarray}
which will prove the identity for Eq.~\eqref{eq:94} when $k = m \times p$. The equality in $=_1$ will be proved by the following Lemma~\ref{lma:15}.
$\hfill \Box$

\begin{lemma}\label{lma:15}
Let $\mathcal{D}_\tau$ be TPSD T-product tensors with $\int_{\Omega} \left\Vert \mathcal{D}_{\tau}^{-p}\right\Vert_\rho d \nu (\tau) < \infty$ for any $p > 0$, then we have
\begin{eqnarray}\label{eq1:lma:15}
\lim\limits_{p \rightarrow 0} \left( \frac{1}{p} \log \int_{\Omega} \frac{ \left\Vert \mathcal{D}_\tau^{-p}\right\Vert_1 }{m \times p} d \nu (\tau)\right) &=& -\frac{1}{m \times p} \int_{\Omega} \log \det( \mathcal{D}_{\tau} ) d \nu (\tau)
\end{eqnarray}
\end{lemma}
\textbf{Proof:}
Because $\int_{\Omega} \left\Vert \mathcal{D}_{\tau}^{-p}\right\Vert_\rho d \nu (\tau) < \infty$, we have that $\mathcal{D}_{\tau}$ are TPD tensors for $\tau$ almost everywhere in $\Omega$. Then, we have 
\begin{eqnarray}
\lim\limits_{p \rightarrow 0} \left( \frac{1}{p}\log \int_{\Omega} \frac{ \left\Vert \mathcal{D}_\tau^{-p}\right\Vert_1}{m \times p} d \nu (\tau) \right) &=_1& \lim\limits_{p \rightarrow 0}\frac{   \int_{\Omega} \frac{ - \sum\limits_{i=1}^{m \times p}  \log \lambda_i(\mathcal{D}_\tau)   }{m \times p} d \nu (\tau)   }{  \int_{\Omega} \frac{ \left\Vert \mathcal{D}_\tau^{-p}\right\Vert_1}{m \times p} d \nu (\tau)      } \nonumber \\
&=& \frac{-1}{m \times p} \int_{\Omega} \sum\limits_{i=1}^{m \times p} \log \lambda_i(\mathcal{D}_\tau)  d \nu (\tau)  
\nonumber \\
&=_2& \frac{-1}{m \times p} \int_{\Omega} \log \det ( \mathcal{D}_\tau ) d \nu (\tau), 
\end{eqnarray}
where $=_1$ is from L'Hopital's rule, and $=_2$ is obtained from $\det$ definition.
$\hfill \Box$

\subsection{T-product Tensor Norm Inequalities by Majorization}\label{sec:T-product Tensor Norm Inequalities by Majorization}


In this section, we will apply derived majorization inequalities for T-product tensors to multivariate T-product tensor norm inequalities which will be used to bound random T-product tensor concentration inequalities in later sections. We will begin to present a Lie-Trotter product formula for tensors. 
\begin{lemma}\label{lma: Lie product formula for tensors}
Let $m \in \mathbb{N}$ and $(\mathcal{L}_k)_{k=1}^{m}$ be a finite sequence of bounded T-product tensors with dimensions $\mathcal{L}_k \in  \mathbb{R}^{m \times m \times p}$, then we have
\begin{eqnarray}
\lim_{n \rightarrow \infty} \left(  \prod_{k=1}^{m} \exp(\frac{  \mathcal{L}_k}{n})\right)^{n}
&=& \exp \left( \sum_{k=1}^{m}  \mathcal{L}_k \right)
\end{eqnarray}
\end{lemma}
\textbf{Proof:}

We will prove the case for $m=2$, and the general value of $m$ can be obtained by mathematical induction. 
Let $\mathcal{L}_1, \mathcal{L}_2$ be bounded tensors act on some Hilbert space. Define $\mathcal{C} \define \exp( (\mathcal{L}_1 + \mathcal{L}_2)/n) $, and $\mathcal{D} \define \exp(\mathcal{L}_1/n) \star \exp(\mathcal{L}_2/n)$. Note we have following estimates for the norm of tensors $\mathcal{C}, \mathcal{D}$: 
\begin{eqnarray}\label{eq0: lma: Lie product formula for tensors}
\left\Vert \mathcal{C} \right\Vert, \left\Vert \mathcal{D} \right\Vert \leq \exp \left( \frac{\left\Vert \mathcal{L}_1 \right\Vert + \left\Vert \mathcal{L}_2 \right\Vert  }{n} \right) =  \left[ \exp \left(  \left\Vert \mathcal{L}_1 \right\Vert + \left\Vert \mathcal{L}_2 \right\Vert  \right) \right]^{1/n}.
\end{eqnarray}

From the Cauchy-Product formula, the tensor $\mathcal{D}$ can be expressed as:
\begin{eqnarray}
\mathcal{D} &=& \exp(\mathcal{L}_1/n) \star \exp(\mathcal{L}_2/n) = \sum_{i = 0}^{\infty} \frac{( \mathcal{L}_1/n)^i}{i !} \star \sum_{j = 0}^{\infty} \frac{( \mathcal{L}_2/n)^j}{j !} \nonumber\\
&=& \sum_{l = 0}^{\infty} n^{-l} \sum_{i=0}^l \frac{\mathcal{L}_1^i}{i!} \star \frac{\mathcal{L}_2^{l-i}}{(l - i)!},
\end{eqnarray}
then we can bound the norm of $\mathcal{C} - \mathcal{D}$ as 
\begin{eqnarray}\label{eq1: lma: Lie product formula for tensors}
\left\Vert \mathcal{C} - \mathcal{D} \right\Vert &=& \left\Vert \sum_{i=0}^{\infty} \frac{\left( [ \mathcal{L}_1 + \mathcal{L}_2]/n \right)^i}{i! }
 - \sum_{l = 0}^{\infty} n^{-l} \sum_{i=0}^l \frac{\mathcal{L}_1^i}{i!} \star \frac{\mathcal{L}_2^{l-i}}{(l - i)!} \right\Vert \nonumber \\
&=& \left\Vert \sum_{i=2}^{\infty} k^{-i} \frac{\left( [ \mathcal{L}_1 + \mathcal{L}_2] \right)^i}{i! }
 - \sum_{m = l}^{\infty} n^{-l} \sum_{i=0}^l \frac{\mathcal{L}_1^i}{i!} \star \frac{\mathcal{L}_2^{l-i}}{(l - i)!} \right\Vert \nonumber \\
& \leq & \frac{1}{k^2}\left[ \exp( \left\Vert \mathcal{L}_1 \right\Vert + \left\Vert \mathcal{L}_2 \right\Vert ) + \sum_{l = 2}^{\infty} n^{-l} \sum_{i=0}^l \frac{\left\Vert \mathcal{L}_1 \right\Vert^i}{i!} \cdot \frac{\left\Vert \mathcal{L}_2 \right\Vert^{l-i}}{(l - i)!} \right] \nonumber \\
& = & \frac{1}{n^2}\left[ \exp \left( \left\Vert \mathcal{L}_1 \right\Vert + \left\Vert \mathcal{L}_2 \right\Vert \right) + \sum_{l = 2}^{\infty} n^{-l} \frac{(  \left\Vert \mathcal{L}_1 \right\Vert + \left\Vert \mathcal{L}_2 \right\Vert )^l}{l!} \right] \nonumber \\
& \leq & \frac{2  \exp \left( \left\Vert \mathcal{L}_1 \right\Vert + \left\Vert \mathcal{L}_2 \right\Vert \right) }{n^2}.
\end{eqnarray}

For the difference between the higher power of $\mathcal{C}$ and $\mathcal{D}$, we can bound them as 
\begin{eqnarray}
\left\Vert \mathcal{C}^n - \mathcal{D}^n \right\Vert &=& \left\Vert \sum_{l=0}^{n-1} \mathcal{C}^m (\mathcal{C} - \mathcal{D})\mathcal{C}^{n-l-1} \right\Vert \nonumber \\
& \leq_1 &  \exp ( \left\Vert \mathcal{L}_1 \right\Vert +  \left\Vert \mathcal{L}_2 \right\Vert) \cdot n \cdot \left\Vert \mathcal{L}_1 - \mathcal{L}_2 \right\Vert,
\end{eqnarray}
where the inequality $\leq_1$ uses the following fact 
\begin{eqnarray}
\left\Vert \mathcal{C} \right\Vert^{l} \left\Vert \mathcal{D} \right\Vert^{n - l - 1} \leq \exp \left( \left\Vert \mathcal{L}_1 \right\Vert +  \left\Vert \mathcal{L}_2 \right\Vert \right)^{\frac{n-1}{n}} \leq 
 \exp\left( \left\Vert \mathcal{L}_1 \right\Vert +  \left\Vert \mathcal{L}_2 \right\Vert \right), 
\end{eqnarray}
based on Eq.~\eqref{eq0: lma: Lie product formula for tensors}. By combining with Eq.~\eqref{eq1: lma: Lie product formula for tensors}, we have the following bound
\begin{eqnarray}
\left\Vert \mathcal{C}^n - \mathcal{D}^n \right\Vert &\leq& \frac{2 \exp \left( 2  \left\Vert \mathcal{L}_1 \right\Vert  +  2  \left\Vert \mathcal{L}_2 \right\Vert \right)}{n}.
\end{eqnarray}
Then this lemma is proved when $n$ goes to infity. $\hfill \Box$

Below, new multivariate norm inequalities for T-product tensors are provided according to previous majorization theorems. 
\begin{theorem}\label{thm:Multivaraite Tensor Norm Inequalities}
Let $\mathcal{C}_i \in \mathbb{R}^{m \times m \times p}$ be TPD tensors, where $1 \leq i \leq n$, $\left\Vert \cdot \right\Vert_{\rho}$ be a unitarily invaraint norm with corresponding gauge function $\rho$. For any continous function $f:(0, \infty) \rightarrow [0, \infty)$ such that $x \rightarrow \log f(e^x)$ is convex on $\mathbb{R}$, we have 
\begin{eqnarray}\label{eq1:thm:Multivaraite Tensor Norm Inequalities}
\left\Vert  f \left( \exp \left( \sum\limits_{i=1}^n \log \mathcal{C}_i\right)   \right)  \right\Vert_{\rho} &\leq& \exp \int_{- \infty}^{\infty} \log \left\Vert f \left( \left\vert \prod\limits_{i=1}^{n}  \mathcal{C}_i^{1 + \iota t} \right\vert\right)\right\Vert_{\rho} \beta_0(t) dt ,
\end{eqnarray}
where $\beta_0(t) = \frac{\pi}{2 (\cosh (\pi t) + 1)}$.

For any continous function $g(0, \infty) \rightarrow [0, \infty)$ such that $x \rightarrow g (e^x)$ is convex on $\mathbb{R}$, we have 
\begin{eqnarray}\label{eq2:thm:Multivaraite Tensor Norm Inequalities}
\left\Vert  g \left( \exp \left( \sum\limits_{i=1}^n \log \mathcal{C}_i\right)   \right)  \right\Vert_{\rho} &\leq& \int_{- \infty}^{\infty} \left\Vert g \left( \left\vert \prod\limits_{i=1}^{n}  \mathcal{C}_i^{1 + \iota t} \right\vert\right)\right\Vert_{\rho} \beta_0(t) dt.
\end{eqnarray}
\end{theorem}
\textbf{Proof:}
From Hirschman interpolation theorem~\cite{sutter2017multivariate} and $\theta \in [0, 1]$, we have 
\begin{eqnarray}\label{eq1:Hirschman interpolation}
\log \left\vert h(\theta) \right\vert \leq \int_{- \infty}^{\infty} \log \left\vert h(\iota t) \right\vert^{1 - \theta} \beta_{1 - \theta}(t) d t + \int_{- \infty}^{\infty} \log \left\vert h(1 +  \iota t) \right\vert^{\theta} \beta_{\theta}(t) d t , 
\end{eqnarray}
where $h(z)$ be uniformly bounded on $S \define \{ z \in \mathbb{C}: 0 \leq \Re(z) \leq 1  \}$ and holomorphic on $S$. The term $ d \beta_{\theta}(t) $ is defined as :
\begin{eqnarray}\label{eq:beta theta t def}
\beta_{\theta}(t) \define \frac{ \sin (\pi \theta)}{ 2 \theta (\cos(\pi t) + \cos (\pi \theta))  }.  
\end{eqnarray}
Let $H(z)$ be a uniformly bounded holomorphic function with values in $\mathbb{C}^{m \times m \times p}$. Fix some $\theta \in [0, 1]$ and let $\mathcal{U}, \mathcal{V} \in \mathbb{C}^{m \times m \times p}$ be normalized tensors such that $\langle \mathcal{U}, \mathcal{H}(\theta) \star \mathcal{V} \rangle = \left\Vert H(\theta) \right\Vert$. If we define $h(z)$ as $h(z) \define \langle \mathcal{U}, \mathcal{H}(z) \star \mathcal{V} \rangle $, we have following bound: $\left\vert h(z) \right\vert \leq \left\Vert H(z) \right\Vert $ for all $z \in S$. From Hirschman interpolation theorem, we then have following interpolation theorem for tensor-valued function: 
\begin{eqnarray}\label{eq2:Hirschman interpolation}
\log \left\Vert H(\theta) \right\Vert \leq \int_{- \infty}^{\infty} \log \left\Vert H(\iota t) \right\Vert^{1 - \theta} \beta_{1 - \theta}(t) dt + \int_{- \infty}^{\infty} \log \left\Vert H(1 +  \iota t) \right\Vert^{\theta}  \beta_{\theta}(t) dt .  
\end{eqnarray}

Let $H(z) = \prod\limits_{i=1}^{n} \mathcal{C}^z_i$. Then the first term in the R.H.S. of Eq.~\eqref{eq2:Hirschman interpolation} is zero since $H(\iota t)$ is a product of unitary tensors. Then we have 
\begin{eqnarray}\label{eq:109}
\log \left\Vert \left\vert \prod\limits_{i=1}^{n} \mathcal{C}_i^{\theta} \right\vert^{\frac{1}{\theta}} \right\Vert \leq \int_{- \infty}^{\infty} \log \left\Vert    \prod\limits_{i=1}^{n} \mathcal{C}_i^{1 + \iota t}         \right\Vert \beta_{\theta}(t) d t .  
\end{eqnarray}

From Lemma~\ref{lma:antisymmetric tensor product properties}, we have following relations:
\begin{eqnarray}\label{eq:111-1}
\left\vert \prod\limits_{i=1}^n \left( \wedge^k \mathcal{C}_i \right)^{\theta}\right\vert^{\frac{1}{\theta}} = \wedge^k \left\vert  \prod\limits_{i=1}^n \mathcal{C}^{\theta}_i  \right\vert^{\frac{1}{\theta}},
\end{eqnarray}
and
\begin{eqnarray}\label{eq:111-2}
\left\vert \prod\limits_{i=1}^n \left( \wedge^k \mathcal{C}_i \right)^{1 + \iota t} \right\vert = \wedge^k \left\vert  \prod\limits_{i=1}^n \mathcal{C}^{1 + \iota t}_i  \right\vert.
\end{eqnarray}
If Eq.~\eqref{eq:109} is applied to $\wedge^k \mathcal{C}_i$ for $1 \leq k \leq r$, we have following log-majorization relation from Eqs.~\eqref{eq:111-1} and~\eqref{eq:111-2}:
\begin{eqnarray}\label{eq:112}
\log \vec{\lambda} \left(  \left\vert \prod\limits_{i=1}^{n} \mathcal{C}_i^{\theta} \right\vert^{\frac{1}{\theta}}   \right) \prec \int_{- \infty}^{\infty} \log \vec{\lambda} \left\vert \prod\limits_{i=1}^{n} \mathcal{C}_i^{1 + \iota t } \right\vert^{\frac{1}{\theta}}   \beta_{\theta}(t) d t.
\end{eqnarray}
Moreover, we have the equality condition in Eq.~\eqref{eq:112} for $k = r$ due to following identies:
\begin{eqnarray}\label{eq:113}
\det \left\vert \prod\limits_{i=1}^n \mathcal{C}_i^{\theta} \right\vert^{\frac{1}{\theta}}
= \det \left\vert \prod\limits_{i=1}^n \mathcal{C}_i^{1 + \iota t } \right\vert= \prod\limits_{i=1}^n \det \mathcal{C}_i. 
\end{eqnarray}

At this stage, we are ready to apply Theorem~\ref{thm:int log average thm 10}  for the log-majorization provided by Eq.~\eqref{eq:112} to get following facts:
\begin{eqnarray}\label{eq:114}
\left\Vert  f \left( \left\vert \prod\limits_{i=1}^{n} \mathcal{C}_i^{\theta} \right\vert^{\frac{1}{\theta}}  \right)  \right\Vert_{\rho} &\leq& \exp \int_{- \infty}^{\infty} \log \left\Vert f \left( \left\vert \prod\limits_{i=1}^{n}  \mathcal{C}_i^{1 + \iota t} \right\vert\right)\right\Vert_{\rho} \beta_{\theta}(t) d t ,
\end{eqnarray}
and
\begin{eqnarray}\label{eq:115}
\left\Vert  g \left( \left\vert \prod\limits_{i=1}^{n} \mathcal{C}_i^{\theta} \right\vert^{\frac{1}{\theta}}  \right)  \right\Vert_{\rho} &\leq& \int_{- \infty}^{\infty} \left\Vert g \left( \left\vert \prod\limits_{i=1}^{n}  \mathcal{C}_i^{1 + \iota t} \right\vert\right)\right\Vert_{\rho}  \beta_{\theta}(t) d t.
\end{eqnarray}
From Lie product formula for tensors given by Lemma~\ref{lma: Lie product formula for tensors},  we have 
\begin{eqnarray}\label{eq:117}
 \left\vert \prod\limits_{i=1}^{n} \mathcal{C}_i^{\theta} \right\vert^{\frac{1}{\theta}}
\rightarrow \exp \left(  \sum\limits_{i=1}^{n} \log \mathcal{C}_i  \right). 
\end{eqnarray}
By setting $\theta \rightarrow 0$ in Eqs.~\eqref{eq:114},~\eqref{eq:115} and using Lie product formula given by Eq.~\eqref{eq:117},  we will get Eqs.~\eqref{eq1:thm:Multivaraite Tensor Norm Inequalities} and~\eqref{eq2:thm:Multivaraite Tensor Norm Inequalities}. 
$\hfill \Box$

\section{Applications of T-product Tensor Norm Inequalities}\label{sec:Applications of T-product Tensor Norm Inequalities}

The purpose of this section is to apply new derived T-product tensor norm inequalities to obtain random symmetric T-product tensor Bernstein bounds. In Section~\ref{sec:Ky Fan k norm Concentration Bound}, Ky Fan $k$-norm inequalities for T-product tensors will be provided and such Ky Fan $k$-norm inequalities will be utilized to establish T-product tensor Bernstein bounds in Section~\ref{sec:T-product Tensor Bernstein Bound} and Section~\ref{sec:Generalized T-product Tensor Bernstein Bound}.

\subsection{Ky Fan $k$-norm Tail Bounds}\label{sec:Ky Fan k norm Concentration Bound}


We will present several lemmas required to prove Ky Fan $k$-norm tail bounds. 
\begin{lemma}\label{lma:T eigenvalue sum rep}
Given a symmetric T-product tensor $\mathcal{C} \in \mathbb{R}^{m \times m \times p}$ which can be expressed as the format shown by Eq.~\eqref{eq:block diagonalized format}, the T-eigenvalues of  $\mathcal{C}$ with respect to the matrix $\mathbf{C}_i$ are denoted as $\lambda_{i, k_i}$, where $1 \leq k_i \leq m$, and we assume that 
$\lambda_{i, 1} \geq \lambda_{i, 2} \geq \cdots \geq \lambda_{i, m}$ (including multiplicities). We have following relation about T-eigenvalues summation representation:
\begin{eqnarray}
\sum\limits_{i=1}^{p} \max\limits_{\mathbf{U}_i \mathbf{U}^{\mathrm{H}}_i = \mathbf{I}_{k_i }} \mathrm{Tr} \mathbf{U}_i \mathbf{C}_i \mathbf{U}^{\mathrm{H}}_i = \sum\limits_{i=1}^{p} \sum\limits_{j=1}^{k_i} \lambda_{i, j}(\mathbf{C}_i),
\end{eqnarray}
and 
\begin{eqnarray}
\sum\limits_{i=1}^{p} \min\limits_{\mathbf{U}_i \mathbf{U}^{\mathrm{H}}_i  = \mathbf{I}_{k_i}} \mathrm{Tr} \mathbf{U}_i \mathbf{C}_i \mathbf{U}^{\mathrm{H}}_i = \sum\limits_{i=1}^{p} \sum\limits_{j=1}^{k_i} \lambda_{i, m - j +1}(\mathbf{C}_i). 
\end{eqnarray}
\end{lemma}
\textbf{Proof:}
From Theorem~\ref{thm:tensor and matrix pd relation}, each matrix $\mathbf{C}_i$ associated to $\mathcal{C}$ based on the format shown by Eq.~\eqref{eq:block diagonalized format} is Hermitian, then the matrix $\mathbf{C}_i$ can be diagonalized as $\mathbf{D}_i$ by the unitary matrix $\mathbf{V}_i$. Without loss of generality, we may assume that $\mathbf{C}_i$ are diagonal matrices. Then, we have
\begin{eqnarray}\label{eq1:lma:T eigenvalue sum rep}
\mathrm{Tr} \mathbf{U}_i \mathbf{D}_i \mathbf{U}^{\mathrm{H}}_i = \sum\limits_{j=1}^{k_i}
 \sum\limits_{l=1}^{m} u_{i_{k, l}}u^{\ast}_{i_{k, l}} \lambda_{i, l} = \sum\limits_{j=1}^{k_i}
 \sum\limits_{l=1}^{m} p_{i_{k, l}}\lambda_{i, l} = [\overbrace{1, \cdots, 1}^{k_i \mbox{terms}}] \mathbf{P}_i
\left[
    \begin{array}{c}
       \lambda_{i,1} \\
       \lambda_{i,2} \\
       \vdots    \\
       \lambda_{i,m}\\
    \end{array}
\right], 
\end{eqnarray}
where the superscript $\ast$ is the operation of a complex conjugate, and $\mathbf{P}_i= (p_{i_{k, l}})$ is a $k_i \times m$ stochastic matrix. From the fact provided by 2.C.1 in~\cite{marshall2011inequalities}, there exists an $(m - k_i) \times m$ matrix $\mathbf{Q}_i$ such that  $ \left[ \begin{array}{c}
       \mathbf{P}_i  \\
       \mathbf{Q}_i \\
    \end{array}
\right]$ is a doubly stochastic matrix. Then, Eq.~\eqref{eq1:lma:T eigenvalue sum rep} can be expressed as 
\begin{eqnarray}
\mathrm{Tr} \mathbf{U}_i \mathbf{D}_i \mathbf{U}^{\mathrm{H}}_i &=&
 [\overbrace{1, \cdots, 1}^{k_i \mbox{terms}}, \overbrace{0, \cdots,0}^{m - k_i \mbox{terms}}] 
\left[ \begin{array}{c}
       \mathbf{P}_i  \\
       \mathbf{Q}_i \\
    \end{array}
\right] \left[
    \begin{array}{c}
       \lambda_{i,1} \\
       \lambda_{i,2} \\
       \vdots    \\
       \lambda_{i,m}\\
    \end{array}
\right].
\end{eqnarray}
Because we have 
\begin{eqnarray}
[ \lambda_{i,1}, \cdots,  \lambda_{i,m}] [ \mathbf{P}^{\mathrm{T}}_i  \mathbf{Q}^{\mathrm{T}}_i] \preceq 
[ \lambda_{i,1}, \cdots,  \lambda_{i,m}],
\end{eqnarray}
then, we can apply the fact 3.H.2 about majorization in~\cite{marshall2011inequalities} to get
\begin{eqnarray}\label{eq2:lma:T eigenvalue sum rep}
\mathrm{Tr} \mathbf{U}_i \mathbf{D}_i \mathbf{U}^{\mathrm{H}}_i \leq  \sum\limits_{j=1}^{k_i} \lambda_{i, j}(\mathbf{C}_i),
\end{eqnarray}
and 
\begin{eqnarray}\label{eq3:lma:T eigenvalue sum rep}
\mathrm{Tr} \mathbf{U}_i \mathbf{D}_i \mathbf{U}^{\mathrm{H}}_i \geq \sum\limits_{j=1}^{k_i} \lambda_{i, m - j +1}(\mathbf{C}_i),
\end{eqnarray}
where $1 \leq k_i \leq m$. We have to note that Eqs~\eqref{eq2:lma:T eigenvalue sum rep} and~\eqref{eq3:lma:T eigenvalue sum rep} are achieved with equalities for $\mathbf{U}_i \mathbf{V}_i = (\mathbf{I}_{k_i}, \mathbf{O})$ and $\mathbf{U}_i \mathbf{V}_i = ( \mathbf{O}, \mathbf{I}_{k_i})$, respectively. This lemma is proved by taking summation with respect to the index $i$ at both sides of Eqs~\eqref{eq2:lma:T eigenvalue sum rep} and~\eqref{eq3:lma:T eigenvalue sum rep}.
$\hfill \Box$

Following lemma will apply Lemma~\ref{lma:T eigenvalue sum rep} to prove majorization relation between T-product tensors summation. 
\begin{lemma}\label{lma:singular values major relation}
Given two symmetric T-product tensors $\mathcal{C}, \mathcal{D} \in \mathbb{R}^{m \times m \times p}$. We have following majorization relation about T-singular values:
\begin{eqnarray}
\sigma(\mathcal{C} +\mathcal{D}) \prec_{w} \sigma(\mathcal{C}) + \sigma(\mathcal{D}).
\end{eqnarray}
\end{lemma}
\textbf{Proof:}
Since we have
\begin{eqnarray}
\sum\limits_{i=1}^{p} \sum\limits_{j=1}^{k_i} \sigma_{i, j}(\mathcal{C} + \mathcal{D}) &=_1& 
\max\limits_{\mathbf{U}_i \mathbf{U}^{\mathrm{H}}_i =_1 \mathbf{I}_{k_i }} \Re \left( \sum\limits_{i=1}^{p}  \mathrm{Tr} \mathbf{U}_i (\mathbf{C}_i + \mathbf{D}_i) \mathbf{U}^{\mathrm{H}}_i \right) \nonumber \\
&\leq&  \max\limits_{\mathbf{U}_i \mathbf{U}^{\mathrm{H}}_i =_1 \mathbf{I}_{k_i }} \Re \left( \sum\limits_{i=1}^{p}  \mathrm{Tr} \mathbf{U}_i \mathbf{C}_i \mathbf{U}^{\mathrm{H}}_i \right) +
\max\limits_{\mathbf{U}_i \mathbf{U}^{\mathrm{H}}_i =_1 \mathbf{I}_{k_i }} \Re \left( \sum\limits_{i=1}^{p}  \mathrm{Tr} \mathbf{U}_i \mathbf{D}_i \mathbf{U}^{\mathrm{H}}_i \right)\nonumber \\
&=_2&  \sum\limits_{i=1}^{p} \sum\limits_{j=1}^{k_i} \sigma_{i, j}(\mathcal{C}) +
\sum\limits_{i=1}^{p} \sum\limits_{j=1}^{k_i} \sigma_{i, j}(\mathcal{D}) 
\end{eqnarray}
where $\Re$ is the operation to take the real part, and the equalities $=_1$ and $=_2$ come from Lemma~\ref{lma:T eigenvalue sum rep}.
$\hfill \Box$

We are ready to introduce the following two lemmas about Ky Fan $k$-norm inequalities for the product of tensors (Lemma~\ref{lma:Ky Fan Inequalities for the prod of tensosrs}) and the summation of tensors (Lemma~\ref{lma:Ky Fan Inequalities for the sum of tensosrs}).

\begin{lemma}\label{lma:Ky Fan Inequalities for the prod of tensosrs}
Let $\mathcal{C}_i \in \mathbb{R}^{m \times m \times p}$ be symmetric T-product tensorsand let $p_i$ be positive real numbers satisfying $\sum\limits_{i=1}^m \frac{1}{p_i} = 1$. Then, we have 
\begin{eqnarray}\label{eq1:lma:Ky Fan Inequalities for the prod of tensosrs}
\left\Vert \left\vert  \prod\limits_{i=1}^{m} \mathcal{C}_i \right\vert^s \right\Vert_{(k)}
\leq  \prod\limits_{i=1}^{m} \left(  \left\Vert \left\vert \mathcal{C}_i \right\vert^{s p_i} \right\Vert_{(k)}    \right)^{\frac{1}{p_i}} 
\leq  \sum\limits_{i=1}^{m} \frac{ \left\Vert \left\vert \mathcal{C}_i \right\vert^{s p_i} \right\Vert_{(k)}      }{p_i}
\end{eqnarray}
where $s \geq 1$ and $k \in \{1,2,\cdots, m \times p  \}$. 
\end{lemma}
\textbf{Proof:}
Since we have 
\begin{eqnarray}
\left\Vert \left\vert  \prod\limits_{i=1}^{m} \mathcal{C}_i \right\vert^s \right\Vert_{(k)} 
= \sum\limits_{j=1}^k \lambda_j \left(   \left\vert  \prod\limits_{i=1}^{m} \mathcal{C}_i \right\vert^s    \right)  = \sum\limits_{j=1}^k \lambda^s_j \left(   \left\vert  \prod\limits_{i=1}^{m} \mathcal{C}_i \right\vert    \right) =
\sum\limits_{j=1}^k \sigma^s_j \left(    \prod\limits_{i=1}^{m} \mathcal{C}_i    \right),
\end{eqnarray}
where we have orders for eigenvalues as $\lambda_1 \geq \lambda_2 \geq \cdots $, and singular values as  $\sigma_1 \geq \sigma_2 \geq \cdots $. 

From Lemma~\ref{lma:antisymmetric tensor product properties}, we have 
\begin{eqnarray}
\left\Vert \left(\prod\limits_{i=1}^{m} \mathcal{C}_i \right)^{\wedge k} \right\Vert 
= \prod\limits_{j=1}^{k}  \sigma_j \left( \prod\limits_{i=1}^{m} \mathcal{C}_i  \right).
\end{eqnarray}
Apply Theorem H.1. in ~\cite{marshall2011inequalities} to each matrix at block diagonal of $\mbox{bcirc}(\mathcal{C}_i) $ by Eq.~\eqref{eq:block diagonalized format}, we will have 
\begin{eqnarray}\label{eq2:lma:Ky Fan Inequalities for the prod of tensosrs}
\sum\limits_{j=1}^k   \sigma^s_j \left( \prod\limits_{i=1}^{m} \mathcal{C}_i  \right)
 \leq \sum\limits_{j=1}^k \left( \prod\limits_{i=1}^{m}   \sigma^s_j (\mathcal{C}_i ) \right).
\end{eqnarray}
Then, we can apply H\"older's inequality to Eq.~\eqref{eq2:lma:Ky Fan Inequalities for the prod of tensosrs} and obtain 
\begin{eqnarray}
\sum\limits_{j=1}^k \left( \prod\limits_{i=1}^{m}   \sigma^s_j (\mathcal{C}_i ) \right) 
&\leq& \prod_{i=1}^{m} \left( \sum\limits_{j=1}^k \sigma^{s p_i}_j (\mathcal{C}_i) \right)^{\frac{1}{p_i}} = \prod_{i=1}^{m} \left( \sum\limits_{j=1}^k \lambda^{s p_i}_j ( \left\vert \mathcal{C}_i \right\vert ) \right)^{\frac{1}{p_i}}  \nonumber \\
&=& \prod_{i=1}^{m} \left( \sum\limits_{j=1}^k \lambda_j ( \left\vert \mathcal{C}_i \right\vert^{s p_i} ) \right)^{\frac{1}{p_i}}= \prod\limits_{i=1}^{m} \left(  \left\Vert \left\vert \mathcal{C}_i \right\vert^{s p_i} \right\Vert_{(k)}    \right)^{\frac{1}{p_i}}
\end{eqnarray}

The second inequality in Eq.~\eqref{eq1:lma:Ky Fan Inequalities for the prod of tensosrs} is obtained by applying Young's inequality to numbers $  \left\Vert \left\vert \mathcal{C}_i \right\vert^{s p_i} \right\Vert_{(k)} $ for $1 \leq i \leq m$. This completes the proof.
$\hfill \Box$

\begin{lemma}\label{lma:Ky Fan Inequalities for the sum of tensosrs}
Let $\mathcal{C}_i \in \mathbb{R}^{m \times m \times p}$ be symmetric T-product tensors, then we have 
\begin{eqnarray}\label{eq1:lma:Ky Fan Inequalities for the sum of tensosrs}
\left\Vert \left\vert  \sum\limits_{i=1}^{m} \mathcal{C}_i \right\vert^s \right\Vert_{(k)}
\leq  m^{s -1} \sum\limits_{i=1}^{m}  \left\Vert \left\vert \mathcal{C}_i \right\vert^{s} \right\Vert_{(k)}     
\end{eqnarray}
where $s \geq 1$ and $k \in \{1,2,\cdots, m \times p\}$. 
\end{lemma}
\textbf{Proof:}
Since we have 
\begin{eqnarray}
\left\Vert \left\vert  \sum\limits_{i=1}^{m} \mathcal{C}_i \right\vert^s \right\Vert_{(k)}
 = \sum\limits_{j=1}^{k} \lambda_j \left( \left\vert  \sum\limits_{i=1}^{m} \mathcal{C}_i \right\vert^s \right) =  \sum\limits_{j=1}^{k} \lambda^s_j \left( \left\vert  \sum\limits_{i=1}^{m} \mathcal{C}_i \right\vert \right) = \sum\limits_{j=1}^{k} \sigma^s_j \left( \sum\limits_{i=1}^{m} \mathcal{C}_i \right). 
\end{eqnarray}
where we have orders for eigenvalues as $\lambda_1 \geq \lambda_2 \geq \cdots $, and singular values as  $\sigma_1 \geq \sigma_2 \geq \cdots $. 

From Lemma~\ref{lma:singular values major relation}, we have
%
\begin{eqnarray}
\sum\limits_{j = 1}^{k} \sigma_j ( \sum\limits_{i=1}^m \mathcal{C}_i ) \leq 
\sum\limits_{j = 1}^{k} \left( \sum\limits_{i=1}^m \sigma_j (\mathcal{C}_i ) \right),
\end{eqnarray}
where $k \in \{1,2,\cdots, m \times p\}$. Then, we have 
\begin{eqnarray}
\sum\limits_{j = 1}^{k} \sigma^s_j ( \sum\limits_{i=1}^m \mathcal{C}_i ) &\leq &
\sum\limits_{j = 1}^{k} \left( \sum\limits_{i=1}^m \sigma_j (\mathcal{C}_i ) \right)^s
\leq m^{s-1} \sum\limits_{j = 1}^{k} \left( \sum\limits_{i=1}^m \sigma^s_j (\mathcal{C}_i ) \right) \nonumber \\
& = & 
m^{s-1} \sum\limits_{j = 1}^{k} \left( \sum\limits_{i=1}^m \sigma^s_j ( \left\vert \mathcal{C}_i \right\vert ) \right) = m^{s-1} \sum\limits_{j = 1}^{k} \left( \sum\limits_{i=1}^m \sigma_j ( \left\vert \mathcal{C}_i \right\vert^s ) \right) \nonumber \\
& = & m^{s -1} \sum\limits_{i=1}^{m}  \left\Vert \left\vert \mathcal{C}_i \right\vert^{s} \right\Vert_{(k)}    
\end{eqnarray}
$\hfill \Box$

Now, we are ready to present our main theorem about Ky Fan $k$-norm probability bound for a function of tensors summation. 
\begin{theorem}\label{thm:Ky Fan norm prob bound for fun of tensors sum}
Consider a sequence $\{ \mathcal{X}_j  \in \mathbb{R}^{m \times m \times p} \}$ of independent, random, symmetric T-product tensors. Let $g$ be a polynomial function with degree $n$ and nonnegative coefficients $a_0, a_1, \cdots, a_n$ raised by power $s \geq 1$, i.e., $g(x) = \left(a_0 + a_1 x  +\cdots + a_n x^n \right)^s$. Suppose following condition is satisfied:
\begin{eqnarray}\label{eq:special cond}
g \left( \exp\left(t \sum\limits_{j=1}^{m} \mathcal{X}_j \right)\right)  \succeq \exp\left(t g \left( \sum\limits_{j=1}^{m} \mathcal{X}_j   \right) \right)~~\mbox{almost surely},
\end{eqnarray}
where $t > 0$. Then, we have 
\begin{eqnarray}\label{eq1:thm:Ky Fan norm prob bound for fun of tensors sum}
\mathrm{Pr} \left( \left\Vert g\left( \sum\limits_{j=1}^{m} \mathcal{X}_j  \right)\right\Vert_{(k)}  \geq \theta \right) &\leq&  
(n+1)^{s-1}\inf\limits_{t, p_j} \exp\left( - \theta t \right) \nonumber \\
&  & \cdot \left(k a_0^s + \sum\limits_{l=1}^{n} \sum\limits_{j=1}^m \frac{  a^{ls}_l  \mathbb{E} \left\Vert \exp\left( p_j  l s t \mathcal{X}_j \right) \right\Vert_{(k)} }{p_j}     \right),
\end{eqnarray}
where $\sum\limits_{j=1}^m \frac{1}{p_j} =1$ and $p_j > 0$. 
\end{theorem}
\textbf{Proof:}
Let $t > 0$ be a parameter to be chosen later. Then
\begin{eqnarray}\label{eq2:thm:Ky Fan norm prob bound for fun of tensors sum}
\mathrm{Pr} \left( \left\Vert g\left( \sum\limits_{j=1}^{m} \mathcal{X}_j  \right)\right\Vert_{(k)} \geq \theta \right)=  \mathrm{Pr} \left(   \left\Vert \exp\left(   t g\left( \sum\limits_{j=1}^{m} \mathcal{X}_j  \right)\right)\right\Vert_{(k)} \geq \exp\left(\theta t \right) \right) \nonumber \\
\leq_1 \exp \left(- \theta t \right) \mathbb{E} \left(   \left\Vert \exp\left(   t g\left( \sum\limits_{j=1}^{m} \mathcal{X}_j  \right)\right)\right\Vert_{(k)} \right) \nonumber \\
\leq_2  \exp \left(- \theta t \right) \mathbb{E} \left(   \left\Vert g \left( \exp\left(   t  \sum\limits_{j=1}^{m} \mathcal{X}_j  \right)\right)\right\Vert_{(k)} \right)
\end{eqnarray}
where $\leq_1$ uses Markov's inequality, $\leq_2$ requires conditions provided by Eq.~\eqref{eq:special cond}.

We can further bound the expectation term in Eq.~\eqref{eq1:thm:Ky Fan norm prob bound for fun of tensors sum} as 
\begin{eqnarray}\label{eq3:thm:Ky Fan norm prob bound for fun of tensors sum}
\mathbb{E} \left(   \left\Vert g \left( \exp\left(   t  \sum\limits_{j=1}^{m} \mathcal{X}_j  \right)\right)\right\Vert_{(k)} \right) &\leq_3 &
\mathbb{E} \int_{- \infty}^{\infty} \left\Vert  g\left(   \left\vert  \prod\limits_{j=1}^m e^{ (1 + \iota \tau)t \mathcal{X}_j }   \right\vert \right)  \right\Vert_{(k)} \beta_0(\tau) d \tau \nonumber \\
& \leq_4 &  (n+1)^{s-1}\left( k a_0^s  \right. \nonumber \\
&  & \left. + \sum\limits_{l=1}^{n}a^{ls}_l  \mathbb{E} 
\int_{- \infty}^{\infty} \left\Vert \left\vert \prod\limits_{j=1}^m e^{ (1 + \iota \tau) t \mathcal{X}_j }    \right\vert^{ls} \right\Vert_{(k)}  \beta_0(\tau)  d \tau     \right),~~~
\end{eqnarray}
where $\leq_3$ from Eq.~\eqref{eq2:thm:Multivaraite Tensor Norm Inequalities} in Theorem~\ref{thm:Multivaraite Tensor Norm Inequalities}, $\leq_4$ is obtained from function $g$ definition and Lemma~\ref{lma:Ky Fan Inequalities for the sum of tensosrs}. Again, the expectation term in Eq.~\eqref{eq3:thm:Ky Fan norm prob bound for fun of tensors sum} can be further bounded by Lemma~\ref{lma:Ky Fan Inequalities for the prod of tensosrs} as 
\begin{eqnarray}\label{eq4:thm:Ky Fan norm prob bound for fun of tensors sum}
\mathbb{E} \int_{- \infty}^{\infty} \left\Vert \left\vert  \prod\limits_{j=1}^m e^{ (1 + \iota \tau)  t\mathcal{X}_j }    \right\vert^{ls} \right\Vert_{(k)}  \beta_0(\tau) d \tau & \leq & \mathbb{E} \int_{- \infty}^{\infty}  \sum\limits_{j=1}^m \frac{ \left\Vert \left\vert e^{t \mathcal{X}_j } \right\vert^{p_j ls} \right\Vert_{(k)}  }{p_j}     \beta_0( \tau) d \tau  \nonumber \\
& = & \sum\limits_{j=1}^m \frac{ \mathbb{E} \left\Vert e^{p_j  l s t \mathcal{X}_j } \right\Vert_{(k)} }{p_j}. 
\end{eqnarray}
Note that the final equality is obtained due to that the integrand is independent of the variable $\tau$ and \\ $\int_{- \infty}^{\infty} \beta_0(\tau) d \tau= 1$. 

Finally, this theorem is established from Eqs.~\eqref{eq2:thm:Ky Fan norm prob bound for fun of tensors sum},~\eqref{eq3:thm:Ky Fan norm prob bound for fun of tensors sum}, and~\eqref{eq4:thm:Ky Fan norm prob bound for fun of tensors sum}. 
$\hfill \Box$

\textbf{Remarks:} The condition provided by Eq.~\eqref{eq:special cond} can be achieved by normalizing tensors $\mathcal{X}_j$ through scaling.


\subsection{T-product Tensor Bernstein Bound}\label{sec:T-product Tensor Bernstein Bound} 


In this section, we will present a tensor Bernstein bound for the maximum and the minimum T-eigenvalue for summation of random symmetric T-product tensors. We will provide the following definition to define a random structure for the T-product tensor $\mathcal{X} \in \mathbb{R}^{m \times m \times p}$.
\begin{definition}\label{def:Random structure of discussed random symmetric T-product tensor}
Random structure for random symmetric T-product tensor $\mathcal{X} \in \mathbb{R}^{m \times m \times p}$
\begin{enumerate}
    \item \label{itm:1}There are $p$ Hermitian matrices with size $m \times m$, denoted as $\mathbf{X}_1, \mathbf{X}_2, \cdots, \mathbf{X}_p$, obtained from Eq.~\eqref{eq:block diagonalized format}. The entries for the matrix $\mathbf{X}_i$ are denoted by $(x_{i_{j,k}})$, where $x_{i_{j,k}}$ is a complex number.  
    \item \label{itm:2} For each $\mathbf{X}_i$, the random variables  $x_{i_{j, j}}$, $\Re x_{i_{j,k}}$ for $j < k$, and $\Im x_{i_{j,k}}$ for $j < k$, are independent. 
    \item \label{itm:3} For each $\mathbf{X}_i$, the random variables  $x_{i_{j, j}}$ follow Gaussian distribution with zero mean and variance as $\frac{1}{m}$.
    \item \label{itm:4} For each $\mathbf{X}_i$, the random variables $\Re x_{i_{j,k}}$ for $j < k$, and $\Im x_{i_{j,k}}$ for $j < k$, follow Gaussian distribution with zero mean and variance as $\frac{1}{2m}$.
\end{enumerate}
\end{definition}


Following lemma is about the expectation of the largest T-eigenvalue of symmetric T-product tensor $\exp (\gamma \mathcal{X})$, where $\gamma$ is a real number. 
\begin{lemma}\label{lma:exp of largest T eigen tensor exponential}
Given a random symmetric T-product tensor $\mathcal{X} \in \mathbb{R}^{m \times m \times p}$ satisfying Definition~\ref{def:Random structure of discussed random symmetric T-product tensor} and any real number $\gamma$, we have
\begin{eqnarray}\label{eq:lma:exp of largest T eigen tensor exponential}
\mathbb{E}\lambda_1 \left( \exp ( \gamma \mathcal{X}  )\right) \leq  \frac{3 m c_1 c_2}{2}   \int^{\infty}_{- \infty}
(y-2)^{1/2} \exp\left[ \gamma y + c_2 m (y-2)^{3/2} \right] dy \define \Psi(m, \gamma, c_1, c_2)
\end{eqnarray}
where $\lambda_1$ is the largest T-eigenvalue, and $c_1, c_2$ are constants related to the bound of cumulative distribution function of the largest eigenvalue of the random Hermitian matrix $\mathbf{X}$. 
\end{lemma}
\textbf{Proof:}
From random structure of discussed random symmetric T-product tensor $\mathcal{X} \in \mathbb{R}^{m \times m \times p}$, all random Hermitian matrix $\mathbf{X}_i$ have same probability density distributions. The maximum T-eigenvalue of $\mathcal{X}$ will be equal to the maximum eigenvalue of $\mathbf{X}_i$, and we use $\mathbf{X}$ to represent any random Hermitian matrix $\mathbf{X}_i$ since they share same distribution. 

From Eq. (2) in~\cite{aubrun2005sharp}, given a $m \times m $ random Hermitian matrix $\mathbf{X}$, we have
\begin{eqnarray}
\mathrm{Pr}\left(\lambda_1 (\mathbf{X}) \leq y\right) \leq 1 - c_1 \exp (-c_2 m (y-2)^{3/2})
\end{eqnarray}
where $c_1, c_2$ are constants related to the bound of cumulative distribution function of the largest eigenvalue of any random Hermitian matrix $\mathbf{X}$. Then, we have
\begin{eqnarray}
\mathbb{E} \lambda_1(\exp (\gamma \mathcal{X})) &=_1& \mathbb{E} \exp (\gamma \lambda_1(\mathcal{X}))
= \int^{\infty}_{- \infty} \exp (\gamma y) d (\mathrm{Pr}\left( \lambda_1( \mathbf{X} ) \leq y   \right)) \nonumber \\
&\leq&   \int^{\infty}_{- \infty} \exp (\gamma y) d \left\{ 1 - c_1 \exp\left[ c_2 m (y-2)^{3/2} \right] \right\} 
\nonumber \\
& = &  \frac{3 m c_1 c_2}{2}   \int^{\infty}_{- \infty}
(y-2)^{1/2} \exp\left[ \gamma y + c_2 m (y-2)^{3/2} \right] dy,
\end{eqnarray}
where $=_1$ comes from the spectral mapping theorem. 
$\hfill \Box$


We are ready to present our theorem about the maximum and the minimum of T-eigenvalue for the summation of random symmetric T-product tensors.
\begin{theorem}[T-product Tensor Bernstein Bound for T-eigenvalue]\label{thm:Generalized Tensor Bernstein Bound}
Consider a sequence $\{ \mathcal{X}_j  \in \mathbb{R}^{m \times m \times p} \}$ of independent, random, symmetric T-product tensors with random structure defined by Definition~\ref{def:Random structure of discussed random symmetric T-product tensor}. Then we have following inequalities: given $\theta_1 > 0$, we have 
\begin{eqnarray}\label{eq1:thm:Generalized Tensor Bernstein Bound}
\mathrm{Pr} \left( \lambda_{\max}\left( \sum\limits_{j=1}^{m} \mathcal{X}_j  \right)  \geq \theta_1 \right) & \leq & 
\inf\limits_{t > 0} \left[ \frac{\exp (- \theta_1 t )}{m} \sum\limits_{j=1}^m \Psi(m, mt, c_1, c_2) \right], 
\end{eqnarray}
and, given $\theta_2 < 0$, we have 
\begin{eqnarray}\label{eq2:thm:Generalized Tensor Bernstein Bound}
\mathrm{Pr} \left( \lambda_{\min}\left( \sum\limits_{j=1}^{m} \mathcal{X}_j  \right)  \leq \theta_2 \right) & \leq & 
\inf\limits_{t > 0} \left[ \frac{\exp (\theta_2 t )}{m} \sum\limits_{j=1}^m \Psi(m, -mt, c_1, c_2) \right]. 
\end{eqnarray}
The $\Psi$ function is defined by Eq.~\eqref{eq:lma:exp of largest T eigen tensor exponential}.
\end{theorem}
\textbf{Proof:}
Since we have 
\begin{eqnarray}
\mathrm{Pr}\left( \lambda_{\max}\left( \sum\limits_{j=1}^{m} \mathcal{X}_j  \right)  \geq \theta_1   \right)
&=_1& \mathrm{Pr}\left( \sigma_{\max}\left( \sum\limits_{j=1}^{m} \mathcal{X}_j  \right)  \geq \theta_1   \right) \nonumber \\
&\leq_2 & \inf\limits_{t, p_j} \exp(- \theta_1 t ) \left( \sum\limits_{j=1}^m \frac{ \mathbb{E}\sigma_{\max} \left( \exp( p_j t \mathcal{X}_j)\right)    }{p_j}      \right)  \nonumber \\
&\leq_3 &  \inf\limits_{t, p_j} \exp(- \theta_1 t ) \left( \sum\limits_{j=1}^m \frac{ \Psi(m, p_j t, c_1, c_2)    }{p_j}      \right)  \nonumber \\
& \leq_4 &  \inf\limits_{t > 0} \exp(- \theta_1 t ) \left( \sum\limits_{j=1}^m \frac{ \Psi(m, m t, c_1, c_2)    }{m}      \right) \nonumber \\
&=& \inf\limits_{t > 0} \left[ \frac{\exp (- \theta_1 t )}{m} \sum\limits_{j=1}^m \Psi(m, mt, c_1, c_2) \right], 
\end{eqnarray}
where $=_1$ comes from that maximum singular value equals to the maximum absolute value of an T-eigenvalue and the maximum and the minimum of T-eigenvalue has same distribution due to the symmetry of random structure given by Definition~\ref{def:Random structure of discussed random symmetric T-product tensor}; the inequality $\leq_2$ comes from Theorem~\ref{thm:Ky Fan norm prob bound for fun of tensors sum} when $g$ is the identity function; the equality $\leq_3$ comes from Lemma~\ref{lma:exp of largest T eigen tensor exponential} and $\sigma_{\max} \left( \exp( p_j t \mathcal{X}_j)\right) = \lambda_{\max} \left( \exp( p_j t \mathcal{X}_j)\right) $ due to TPD of $\exp( p_j t \mathcal{X}_j ) $; the inequality $\leq_4$ is obtained by selecting $p_j = m$. Therefore, we have Eq.~\eqref{eq1:thm:Generalized Tensor Bernstein Bound}.

For the minimum T-eigenvalue, we also have
\begin{eqnarray}
\mathrm{Pr}\left( \lambda_{\min}\left( \sum\limits_{j=1}^{m} \mathcal{X}_j  \right)  \leq \theta_2   \right)
&=_1&  \mathrm{Pr}\left( \lambda_{\max}\left( \sum\limits_{j=1}^{m}- \mathcal{X}_j  \right)  \geq -\theta_2  \right) \nonumber \\
&=_2&\mathrm{Pr}\left( \sigma_{\max}\left( \sum\limits_{j=1}^{m}- \mathcal{X}_j  \right)  \geq -\theta_2  \right) \nonumber \\
&\leq_3& \inf\limits_{t, p_j} \exp( \theta_2 t ) \left( \sum\limits_{j=1}^m \frac{ \mathbb{E}\sigma_{\max} \left( \exp( - p_j t \mathcal{X}_j)\right)    }{p_j}      \right)  \nonumber \\
& \leq_4 &  \inf\limits_{t, p_j} \exp(\theta_2 t ) \left( \sum\limits_{j=1}^m \frac{ \Psi(m, -p_j t, c_1, c_2)    }{p_j}      \right)  \nonumber \\
& \leq_5 & \inf\limits_{t > 0} \left[ \frac{\exp (\theta_2 t )}{m} \sum\limits_{j=1}^m \Psi(m, -mt, c_1, c_2) \right],
\end{eqnarray}
where $=_1$ comes from Theorem~\ref{thm:Courant-Fischer T-product}; $=_2$ is true since the maximum singular value equals to the maximum absolute value of an T-eigenvalue and the maximum and the minimum of T-eigenvalue has same distribution due to the symmetry of random structure given by Definition~\ref{def:Random structure of discussed random symmetric T-product tensor}; the inequality $\leq_3$ comes from Theorem~\ref{thm:Ky Fan norm prob bound for fun of tensors sum} again when $g$ is an identity map; the equality $\leq_4$ comes from Lemma~\ref{lma:exp of largest T eigen tensor exponential} and $\sigma_{\max} \left( \exp( p_j t \mathcal{X}_j)\right) = \lambda_{\max} \left( \exp( p_j t \mathcal{X}_j)\right) $ due to TPD of $\exp( p_j t \mathcal{X}_j)$; the inequality $\leq_5$ is obtained by selecting $p_j = m$. Hence, we have Eq.~\eqref{eq2:thm:Generalized Tensor Bernstein Bound}.
$\hfill \Box$

\subsection{Generalized T-product Tensor Bernstein Bound}\label{sec:Generalized T-product Tensor Bernstein Bound}

In this section, we will present a generalized tensor Bernstein bound for Ky Fan $k$-norm, and we will begin with a lemma to bound exponential of a random T-product tensor.
\begin{lemma}\label{lem:Subexponential Bernstein mgf}
Suppose that $\mathcal{X} \in \mathbb{R}^{m \times m \times p}$ is a random symmetric T-product tensor that satisfies
\begin{eqnarray}\label{eq1:lem:Subexponential Bernstein mgf}
\mathcal{X}^p  \preceq \frac{p! \mathcal{A}^2}{2} 
\mbox{~~ almost surely for $p=2,3,4,\cdots$,} 
\end{eqnarray}
where $\mathcal{A}$ is a fixed TPD tensor. Then, we have
\begin{eqnarray}
 e^{t \mathcal{X}} \preceq \mathcal{I} + t \mathcal{X} + \frac{t^2 \mathcal{A}^2}{2 ( 1 - t) } 
\mbox{~~ almost surely,} 
\end{eqnarray}
where $0 < t < 1$.
\end{lemma}
\textbf{Proof:}
From Tayler series of the tensor exponential expansion, we have 
\begin{eqnarray}
e^{t \mathcal{X}} &=& \mathcal{I} + t \mathcal{X} + \sum\limits_{p=2}^{\infty} \frac{t^p (\mathcal{X}^p)}{p!} \preceq \mathcal{I}+  t \mathcal{X} +  \sum\limits_{p=2}^{\infty} \frac{t^p \mathcal{A}^2}{2} = \mathcal{I}+  t \mathcal{X}  + \frac{t^2 \mathcal{A}^2}{2 (1-t)}.
\end{eqnarray}
Therefore, this Lemma is proved. 
$\hfill \Box$

\begin{lemma}\label{lma:exp of largest T sv tensor exponential}
Given a random symmetric T-product tensor $\mathcal{X} \in \mathbb{R}^{m \times m \times p}$ satisfying Definition~\ref{def:Random structure of discussed random symmetric T-product tensor}, we have
\begin{eqnarray}
\mathbb{E}\sigma_1 \left( \mathcal{X}  \right) \leq   \int^{\infty}_{- 2 }
 d_1 \exp( - d_2 m z^{3/2})dz \define \Phi(m, d_1, d_2).
\end{eqnarray}
where $\sigma_1$ is the largest T-singular value, and $d_1, d_2$ are constants related to the upper bound of the largest eigenvalue of the random Hermitian matrix $\mathbf{X}$. 
\end{lemma}
\textbf{Proof:}
From random structure of discussed random symmetric T-product tensor $\mathcal{X} \in \mathbb{R}^{m \times m \times p}$, all random Hermitian matrix $\mathbf{X}_i$ have same probability density distributions. The maximum T-singular value of $\mathcal{X}$ will be equal to the maximum singular value of $\mathbf{X}_i$, and we use $\mathbf{X}$ to represent any random Hermitian matrix $\mathbf{X}_i$ since they share the same distribution.

From Eq. (8) in~\cite{aubrun2005sharp}, given a $m \times m $ random Hermitian matrix $\mathbf{X}$, we have
\begin{eqnarray}
\mathrm{Pr}\left(\sigma_1 (\mathbf{X}) > y\right) \leq d_1 \exp (-d_2 m (y-2)^{3/2})
\end{eqnarray}
where $d_1, d_2$ are constants related to the upper bound of the largest (or smallest) eigenvalue of any random Hermitian matrix $\mathbf{X}$. Then, we have
\begin{eqnarray}
\mathbb{E} \sigma_1 (\mathcal{X}) &=& \int^{\infty}_0 \mathrm{Pr}\left( \sigma_1 (\mathcal{X}) > y \right) d y \nonumber \\
&\leq &\int^{\infty}_0  d_1 \exp \left[ -d_2 m ( y-2 )^{3/2}  \right] d y \nonumber \\
&=& \int^{\infty}_{- 2 }
 d_1 \exp( - d_2 m z^{3/2})dz 
\end{eqnarray}
$\hfill \Box$

Following lemma is about Ky Fan $k$-norm bound for the exponential of a random T-product tensor with subexponential constraints. 
\begin{lemma}\label{lma:Ky Fan k norm bound for random tensors subexp}
Given a symmetric random T-product tensor $\mathcal{X} \in \mathbb{R}^{m \times m \times p}$ with random structure defined by Definition~\ref{def:Random structure of discussed random symmetric T-product tensor} and 
\begin{eqnarray}\label{eq1:lma:Ky Fan k norm bound for random tensors subexp}
\mathcal{X}^p  \preceq \frac{p! \mathcal{A}^2}{2} 
\mbox{~~ almost surely for $p=2,3,4,\cdots$,} 
\end{eqnarray}
where $\mathcal{A}$ is a TPD T-product tensor. Then, we have following bound about the expectation value of Ky Fan $k$-norm for the random T-product tensor $\exp( \theta \mathcal{X} )$
\begin{eqnarray}\label{eq2:lma:Ky Fan k norm bound for random tensors  subexp}
\mathbb{E} \left\Vert \exp( \theta \mathcal{X} ) \right\Vert_{(k)} &\leq&  k \left[1 + \theta \Phi(m, d_1, d_2)  +  \frac{\theta^2}{ 2 (1 - \theta)} \sigma_1 \left( \mathcal{A}^2\right)  \right]. 
\end{eqnarray}
\end{lemma}
\textbf{Proof:}
From Lemma~\ref{lem:Subexponential Bernstein mgf}, we have 
\begin{eqnarray}
\mathbb{E} \left\Vert \exp( \theta \mathcal{X} ) \right\Vert_{(k)} &=& 
\sum\limits_{l=1}^{k} \mathbb{E} \sigma_l \left( \exp(\theta \mathcal{X})\right) \nonumber \\
&\leq &  \sum\limits_{l=1}^{k} \mathbb{E} \sigma_l \left( \mathcal{I} + \theta \mathcal{X} + 
 \frac{\theta^2 \mathcal{A}^2}{ 2 (1 - \theta)} \right) \leq k \mathbb{E} \sigma_1 \left( \mathcal{I} + \theta \mathcal{X} + 
 \frac{\theta^2 \mathcal{A}^2}{ 2 (1 - \theta)} \right) 
\end{eqnarray}
where $\sigma_l (\cdot )$ is the $l$-th largest T-singular value. 

From Lemma~\ref{lma:singular values major relation}, we have $\sigma_1(\mathcal{A} + \mathcal{B}) \leq \sigma_1(\mathcal{A}) + \sigma_1(\mathcal{B})$ for two symmetric T-product tensors $\mathcal{A}$ and $\mathcal{B}$. Then, we can bound $ \mathbb{E} \sigma_1 \left(\mathcal{I} + \theta \mathcal{X} + 
 \frac{\theta^2 \mathcal{A}^2}{ 2 (1 - \theta)} \right) $ as 
\begin{eqnarray}\label{eq1:lma:Ky Fan k norm bound for random tensors subexp}
 \mathbb{E} \sigma_1 \left(\mathcal{I} + \theta \mathcal{X} + 
 \frac{\theta^2 \mathcal{A}^2}{ 2 (1 - \theta)} \right)& \leq &1 + \theta \mathbb{E}   \sigma_1(\mathcal{X}) +  \frac{\theta^2}{ 2 (1 - \theta)} \sigma_1 \left( \mathcal{A}^2\right)  \nonumber \\
&\leq& 1 + \theta \mathbb{E}  \Phi(m, d_1, d_2 )+  \frac{\theta^2}{ 2 (1 - \theta)} \sigma_1 \left( \mathcal{A}^2\right)
\end{eqnarray}
where we use $\Phi(m, d_1, d_2 )$ from Lemma~\ref{lma:exp of largest T sv tensor exponential} to bound $\mathbb{E}   \sigma_1(\mathcal{X}) $ in the last inequality. This Lemma is proved by multiplying $k$ at Eq.~\eqref{eq1:lma:Ky Fan k norm bound for random tensors subexp}
$\hfill \Box$

We are ready to present our main theorem about the generalized T-product tensor Bernstein bound.


\thmGeneralizedTensorBernsteinBound*
\textbf{Proof:}
Since we have
\begin{eqnarray}
\mathrm{Pr} \left( \left\Vert g\left( \sum\limits_{j=1}^{m} \mathcal{X}_j  \right)\right\Vert_{(k)}  \geq \theta \right)  \leq_1  (n+1)^{s-1}\inf\limits_{t, p_j} e^{- \theta t } \left(k a_0^s + \sum\limits_{l=1}^{n} \sum\limits_{j=1}^m \frac{  a^{ls}_l  \mathbb{E} \left\Vert \exp\left( p_j  l s t \mathcal{X}_j \right) \right\Vert_{(k)} }{p_j}     \right)  \nonumber \\
  \leq_2   (n+1)^{s-1}\inf\limits_{t, p_j} e^{- \theta t }\left\{  k a_0^s + \sum\limits_{l=1}^{n} \sum\limits_{j=1}^m \frac{  a^{ls}_l k \left[1 + p_j lst \Phi(m, d_1, d_2) + 
\frac{(p_j lst)^2      \sigma_1 (\mathcal{A}^2)}{2 (1 - p_j lst)} \right]    }{p_j}     \right\}  \nonumber \\
 \leq_3    (n+1)^{s-1} \inf\limits_{t > 0} e^{- \theta t }k  \left\{a_0^s + \sum\limits_{l=1}^{n} a_l^{l s}\left[1 + mlst \Phi(m, d_1, d_2) + 
\frac{(mlst)^2   \sigma_1(\mathcal{A}^2) }{2(1 - mlst)  },~~
   \right] \right\}.
\end{eqnarray}
where the inequality $\leq_1$ comes from Theorem~\ref{thm:Ky Fan norm prob bound for fun of tensors sum}; 
the inequality $\leq_2$ comes from Lemma~\ref{lma:Ky Fan k norm bound for random tensors subexp}; the inequality $\leq_3$ is obtained by setting $p_j = m$. 
$\hfill \Box$



\section{Conclusions}\label{sec:Conclusions}

This work extend previous work in~\cite{chang2021TProdII} by making following generalizations via majorization techniques:  (1) besides bounds related to extreme values of T-eigenvalues, this works considers more general unitarily invariant norm for T-product tensors; (2) this work derives new bounds for any polynomial function raised by any power greater or equal than one for the summation of random symmetric T-product tensors. We also establish the Courant-Fischer min-max theorem for T-product tensors and marjoization relation for T-singular values which are by-products of our procedure to prove the generalized random T-product  Bernstein bounds. Possible future work about this research is to consider tail bounds behaviors for the summation of random symmetric T-product tensors equipped with other random structures different from random structure provided by Definition~\ref{def:Random structure of discussed random symmetric T-product tensor}.


\bibliographystyle{IEEETran}
\bibliography{TProd_Major_Ind_Bib}

\end{document}